\documentclass[11pt]{article}

\usepackage[margin=1in]{geometry}
\usepackage[T1]{fontenc}
\usepackage[utf8]{inputenc}
\usepackage{lmodern}
\usepackage{microtype}
\usepackage{amsmath,amssymb,amsthm}
\usepackage{enumitem}
\usepackage{booktabs}
\usepackage[pagebackref=true]{hyperref}
\usepackage{xcolor}

\hypersetup{colorlinks=true,linkcolor=blue,citecolor=blue,urlcolor=blue,
  pageanchor=true,hypertexnames=true,plainpages=false}

\newtheorem{theorem}{Theorem}[section]

\newtheorem{remark}[theorem]{Remark}

\newtheorem{corollary}[theorem]{Corollary}

\newcommand{\SoneS}{\ensuremath{\mathrm{S1S}}}
\newcommand{\StwoS}{\ensuremath{\mathrm{S2S}}}

\newcommand{\Pclass}{\ensuremath{\mathsf{P}}}
\newcommand{\NPclass}{\ensuremath{\mathsf{NP}}}
\newcommand{\coNPclass}{\ensuremath{\mathsf{coNP}}}

\newcommand{\BPP}{\ensuremath{\mathsf{BPP}}}
\newcommand{\PP}{\ensuremath{\mathsf{PP}}}

\newcommand{\PSPACE}{\ensuremath{\mathsf{PSPACE}}}
\newcommand{\IP}{\ensuremath{\mathsf{IP}}}

\renewcommand*{\backref}[1]{}
\renewcommand*{\backrefalt}[4]{%
  \ifcase #1 %
  \or (Cited on page~#2.)%
  \else (Cited on pages~#2.)%
  \fi}

\title{Michael Rabin's fingerprints across the theory of computation\footnote{(Random) fingerprints were among Rabin's favorite notions.}}
\author{Avi Wigderson\\ Institute for Advanced Study, Princeton\footnote{\tt avi@ias.edu}}
\date{}

\begin{document}
\maketitle

\begin{abstract}
Michael O. Rabin was a major founder and leader of theoretical computer science. He died on April 14, 2026, at the age of 94. This paper aims to survey the wide scope and lasting impact of his research.
\end{abstract}

\newpage
\tableofcontents
\clearpage

\section{Introduction}
\label{sec:introduction}



The theory of computation, which forms the mathematical underpinnings of computer science, is one of the most important and dynamic disciplines. Since the formalization of the notion of computation by Church, Turing, and Post, the field has constantly expanded its goals and achievements, trying to understand computation, its properties, and its applications in numerous settings and forms. Its models and results have led to remarkable technological revolutions, and have illuminated mathematics and the sciences. Friendly introductions to the field, its diverse parts, and their connections are Harel's \emph{Algorithmics} \cite{Harel1987} (1987) and Wigderson's \emph{Mathematics and Computation}~\cite{Wigderson2019} (2019). Excellent textbooks include Arora and Barak's \emph{Computational Complexity: A Modern Approach} \cite{AroraBarak2009}, Goldreich's \emph{Computational Complexity} \cite{Goldreich2008}, Moore and Mertens' \emph{The Nature of Computation} \cite{MooreMertens2011}, Papadimitriou's \emph{Computational Complexity} \cite{Papadimitriou1994}, and Sipser's \emph{Introduction to the Theory of Computation} \cite{Sipser2013}.


Michael Rabin was one of the giants of this field; in many ways, the evolution of his interests reflected the evolution of the field and shaped many aspects of its development for decades. Rabin described his research path in interviews with Dennis Shasha and David Harel \cite{ShashaRabin2010,RabinHarel2015}, starting with his first encounter with computation. Rabin was a mathematics student at the Hebrew University, focusing on algebra and logic. During his M.Sc.\ studies at the Hebrew University under Jacob Levitzki, he read Kleene's book \cite{Kleene1952} \emph{Introduction to metamathematics}, which described the formative work on computation discussed above. This had a major impact on him. Rabin was also a keen observer of the world around him, and the flurry of activity devoted to building electronic computers compounded that impact. He concluded that the study of computation would lead to new and exciting mathematics, and set out to make good on this prediction. His first step was to go to the heart of this activity and work on his Ph.D. under Church at Princeton, which is the starting point of this story.

The breadth, depth, and impact of Rabin's work are staggering. The list of areas in which his work was foundational is imposing; it includes computability theory, verification and model checking, automata and formal language theory, computational complexity theory, cryptography and distributed computation.

Papers can be foundational in many ways. Rabin's papers on computational complexity theory, probabilistic algorithms, and computational algebra created new fields of study, with ample motivation, basic definitions, long-term goals, and initial results. His paper on finite automata with Scott, which earned them the Turing Award, gave the community a template for what it means to study a computational model in depth. His papers in cryptography introduced novel primitives, assumptions, and models. His papers on verification and distributed computing introduced new techniques for solving fundamental open problems in those fields.

Rabin's distinctive approach reveals itself in his writings. His intellectual pursuit of computational phenomena seeks elegant, simple concepts that capture fundamental issues of lasting significance. His initial works on exploring the power of \emph{nondeterminism} and \emph{randomness} in computation inspire us to this day, and these investigations underlie some of the greatest achievements and open problems of our field. His papers exhibit great vision and powers of abstraction. But despite being a theoretician, Rabin also cared about practice and implementation, and often gave numerical examples of input sizes and runtimes in his papers.


This survey aims to review some of Rabin's central contributions across many subfields of the theory. It is organized in sections around these subfields, roughly in chronological order. In each, I will try to describe some of Rabin's work and impact within its historical context. Many of Rabin's ideas became so ingrained in the field's culture and in undergraduate education that it may sometimes be hard to appreciate that they had to be created, and how novel they were at the time of discovery. To illustrate this, I have tried to supply some historical background. The appendix contains what I hope is a nearly complete list of Rabin's writings, including many technical reports that were never published otherwise.

Some disclaimers are in order. First, this survey is naturally biased by my knowledge, interests, and taste. Second, I focused on readability and accessibility to a general mathematical audience. Thus, in trying to convey the essence of complex frameworks and results, I often oversimplified, neglecting to provide many fine details; I hope this is compensated by the many references to the original texts where the reader can find precise definitions and theorems. These limitations also mean that the survey cannot be a comprehensive historical one, and so does not cover many related works; I hope this is compensated by the books and survey papers I suggested in each section as supplementary reading. In the same vein, there is a large diversity of topics here, and in each there are basic terms that may not be familiar to outsiders; I could not define them all for obvious reasons, but Wikipedia or AI would easily remedy this for the reader. Finally, the history of ideas is fascinating but sometimes also frustrating to explore; it is quite common that ideas arose independently by several people, who may be credited differently in subsequent papers and surveys. I have tried to be as accurate and inclusive as possible in the few relevant cases of ambiguity. 


Let me add a brief personal note. I was fortunate to be a colleague of Rabin at the computer science department of the Hebrew University and had an office across the corridor from him for over a decade. One strong impression I had concerns his enthusiasm for research: a few times Rabin burst into my office to tell me about a new idea he had just had, exhibiting the same excitement as a young graduate student with a first discovery. This attitude toward our field and toward research in general was clearly carried over to the classroom, influencing generations of students at the Hebrew University and at Harvard. I feel that this legacy beautifully complements the nature and impact of his work.



\subsection{Preamble: The dawn of computation}\label{Dawn}

We start with a very short account of the field before Rabin entered it. Excellent historical accounts (with different foci) are the books of Chabert~\cite{Chabert1999}, Davis~\cite{Davis2000} and Kleene~\cite{Kleene1952}.

Like many technological revolutions, the computer revolution had important theoretical foundations. Very early writings on general-purpose ``reasoning'' go back to Leibniz \cite{Leibniz1666}, and writings on actual general-purpose ``programming'' go back to Lovelace's work \cite{MenabreaLovelace1843} with Babbage on his \emph{Analytical Engine}. But the most direct impetus lay in the crisis in the foundations of mathematics at the turn of the twentieth century. As is well known, Hilbert, perhaps the foremost mathematician of his age, had a dream. His ``formalization of mathematics'' project was based on two beliefs. First, every true mathematical fact can be proved in an appropriate axiomatic proof system. Second, whether such a proof exists can be determined automatically by a finite ``mechanical procedure'' \cite{Hilbert1902Problems,HilbertAckermann1928}.



The first belief was shattered by G\"odel's incompleteness theorem \cite{Godel1931}, which showed that there are true facts about the integers that cannot be proved in Peano Arithmetic. More generally, he showed that every sufficiently strong logical system has statements such that neither they nor their negations follow from the axioms.


The second belief, also known as the \emph{Entscheidungsproblem}, was shattered independently by Church \cite{Church1936,Church1936Entscheidungsproblem}, using the lambda calculus, and Turing \cite{Turing1936}, using his model of machines. Turing's model, which Church was quick to call the ``Turing machine'' in a 1937 review \cite{Church1937TuringReview}, mathematically captures Hilbert's vague notion of a ``mechanical procedure''---an essential definition for a mathematical proof that such a procedure does not exist for some task. We will call a function \emph{algorithmically solvable} (and will also use the terms \emph{decidable} or \emph{computable}) if it can be computed by a Turing machine that halts on every input. That this notion of computability is universal---namely, that Turing machines can simulate \emph{any} other form of computation---is known as the Church--Turing thesis, and has stood unchallenged to this day.


To refute Hilbert's Entscheidungsproblem, Turing showed that the \emph{halting problem}---the problem that takes as input a description of a Turing machine and asks whether it halts or runs forever, say on an empty input---is unsolvable. Again, we will alternatively use the terms \emph{undecidable} or \emph{uncomputable} below. It is worth pointing out that the halting problem, proved undecidable by Turing, is really a central problem of the nascent field of computer science; every programmer needs to know whether a program is buggy. After Turing, one can ask:

\begin{quote}
\emph{Are there central problems in mathematics that are algorithmically unsolvable?}
\end{quote}


Suddenly, after Turing's paper, this question became amenable to mathematical analysis. This possibility had, and still has, a major impact on mathematics. Mathematicians have always been interested in designing algorithms for basic problems. Early examples, all developed millennia ago, include Chinese and Egyptian solutions to linear equations, Euclid's algorithm for greatest common divisors in the \emph{Elements}, and many constructive proofs in plane geometry. Indeed, the very word \emph{algorithm} derives from a Latinized form of the name of the ninth-century mathematician al-Khwarizmi, whose works described methods for solving arithmetic and algebraic problems. But what about problems that had evaded algorithms? With the mathematical formulation of computation by Turing machines, the game was now afoot to prove such unsolvability results in all areas of mathematics and to discern the boundary between the decidable and the undecidable.

\section{Between decidability and undecidability: algebra}
\label{sec:computability_algebra}

This section focuses on two papers of Rabin, both grand in scope, that establish a solid foundation connecting computation and algebra. They exhibit the value of the computational lens even for pure mathematical problems involving no computation. The first~\cite{Rabin1958}, from his Ph.D. thesis, contains mainly negative results: basic problems in group theory are proved undecidable. The second~\cite{Rabin1960ComputableAlgebra}, presented soon afterward and published in 1960, contains mainly positive results. It is a foundational contribution to the emerging theory of computable algebra, with algorithms and structural theorems for many basic algebraic constructions.


\subsection{Undecidability of problems on finitely presented groups}

Extensions of undecidability to mathematics followed quickly on the heels of Turing's result on the halting problem. With the existence of such a natural undecidable problem, to prove that a favorite problem is undecidable one ``only'' has to reduce the halting problem to it. Suppose the problem is whether an object in a collection of mathematical structures has a particular property. We design an algorithm that translates any given Turing-machine description into a description of such an object, which has the desired property if and only if the given machine halts\footnote{More general reductions are possible. For example, it suffices to show that an oracle deciding the target problem could be used by a finite algorithm to decide the halting problem.}. In other words, we must show that the given problem for such structures can encode general computation.


One of the first mathematical areas to which this idea was applied was group theory. A natural description of a group is by a \emph{presentation}: a set of generators together with a set of relations, namely pairs of words in the generators that are equal as group elements. We assume throughout that presentations are finite, and fix such a presentation. The most basic problem in this setting, and for almost any other algebraic structure, is the \emph{word problem}: given two words in the generators, are they equal as group elements? Algorithms for various special classes of groups were known, most famously Dehn's algorithm for surface groups~\cite{Dehn1912}, but it remained open whether the general problem was decidable.

The first undecidability results came for semigroups, which need not have inverses or identity elements. Their defining relations allow one to replace locally one subword by another, and therefore look similar to the way a Turing-machine configuration changes after one step. Turing machines, however, do not provide the reversibility implicit in equations between words. This hurdle was overcome independently in 1947 by Markov~\cite{Markov1947} and Post~\cite{Post1947}, who proved the word problem for finitely presented semigroups is undecidable. Turing subsequently strengthened the result to the more restricted class of semigroups with cancellation~\cite{Turing1950Semigroups}.


Groups, having inverses and an identity, are more rigid and thus much harder to use for encoding general computation. Nevertheless, Novikov~\cite{Novikov1955} and Boone~\cite{Boone1957,Boone1959} overcame these obstacles, and proved (via very different constructions) that the word problem for finitely presented groups is undecidable.



Rabin's \emph{Annals} paper gives a sweeping generalization of this result. To explain in what sense, let us return to Turing machines. Recall that Turing proved that a basic property of a program, namely whether it halts on a given input, is undecidable. Rice~\cite{Rice1953} extends this undecidability result to every \emph{semantic} property of a program. Here, ``semantic'' means a property of the partial function computed by a program, and ``nontrivial'' simply means that the property holds for some such functions and fails for others. Informally, Rice's theorem means that nothing interesting about the function computed by an arbitrary program can be decided from the program itself.



By analogy, Rabin~\cite{Rabin1958} proves that a vast range of nontrivial properties of groups given by finite presentations are undecidable\footnote{For experts; Rabin points out that unlike most previous undecidability results, which were combinatorial in nature, his proof uses algebraic reductions, specifically employing direct product with amalgamation.}. Adian obtained closely related results independently~\cite{Adian1955}, and the general result is now called the Adian--Rabin theorem.

\begin{theorem}[Adian--Rabin~\cite{Rabin1958,Adian1955}]\label{thm:adian-rabin}
Every Markov property of finitely presented groups is undecidable.
\end{theorem}

A Markov property\footnote{This name was coined soon after by Boone, as both Adian and Rabin follow Markov's analogous conditions for semigroups.} is one that is \emph{nontrivial and semantic}, with similar meanings as above. Here, the word ``semantic'' is easy to explain: the property should depend only on the group defined by the presentation---not on the presentation itself. In other words, it should be invariant under group isomorphism. The word ``nontrivial'' is subtler and a bit less trivial to explain: as before, we demand that one group has the property, and another does not, but moreover that the latter is not a subgroup of a group that does have the property\footnote{This ``heredity'' condition does not present a problem in many applications. To see why, observe that if the word problem (from which the eventual undecidability reduction is made) is solvable for a presentation of some group $G$, then it is solvable for any presentation of any subgroup $H$ of it. This of course requires that the generators of $H$ can be effectively written as words in the generators of $G$, something which can be done in many interesting cases.}. 


%

Rabin, and also Adian, demonstrate the generality of the theorem through many natural examples. In particular, the following group properties are undecidable from a finite presentation: being trivial, finite, cyclic, commutative, solvable, free, and many more. Their techniques also prove other important undecidability results about groups, including the undecidability of testing whether two finite presentations define isomorphic groups. Rabin then draws a direct, purely group-theoretic consequence: there is no complete system of isomorphism invariants that are \emph{effectively computable}. This strikes a serious blow against natural programs that sought to classify all finitely presented groups by a uniform list of effective algebraic invariants\footnote{I call this corollary ``purely group-theoretic'' because such classification programs existed well before computation; the theory of computation offered both a formal definition of what such a classification should look like (e.g., that its invariants are computable) and a way to show that it is hopeless in general.}.

The Adian--Rabin theorem made undecidability a general obstruction to classification of groups given finite presentations, and is widely used in combinatorial group theory as a guide to obtaining algorithms and characterizations for various classes of groups. This is one of the earliest examples of this impact of computation on the practice and discourse of mathematics, and it spread quickly to other fields of mathematics. We note that more recently complexity-theoretic notions such as \NPclass{}-completeness play a similar guiding role in the search for tractable algorithms and characterizations.



\subsection{Computable algebra}

Algebra has been an algorithmic field from its inception. Some early examples like the GCD algorithm are mentioned in subsection~\ref{Dawn} and numerous others were added over the centuries for constructing algebraic objects. In these classical settings, however, the algebraic structures usually had explicit descriptions (e.g. integers for GCD) that served as the inputs and outputs of the algorithms.

A major new stage arrived with the 1956 paper of Fr\"ohlich and Shepherdson~\cite{FrohlichShepherdson1956} and Rabin's 1960 paper~\cite{Rabin1960ComputableAlgebra}\footnote{An early version of Rabin's work is his July 1957 Cornell Summer Institute notes~\cite{Rabin1957ComputableAlgebraicSystems}.}. Both initiated and developed the theory connecting computability and algebra, differing somewhat in scope and purpose. Fr\"ohlich and Shepherdson were motivated principally by effective procedures in field theory and developed a theory of \emph{explicit fields}. Rabin aimed at a general integration of computation into algebra and developed a computational theory of abstract algebraic systems, including groups, rings, and fields. Their basic notions are closely related (though not necessarily in a computable way). Rabin explained the reasons for his choices, which, as we shall see, can lead to different results. We begin with common aspects of the two papers and then discuss contributions unique to Rabin's.

The basic idea in both approaches is simple once a formal notion of computability became available\footnote{Indeed, it follows earlier attempts at such a formalization by Kronecker and van der Waerden, who worked before it was available.}. A countable abstract algebraic system is equipped with an effective presentation: we index the elements of the system by integers (or binary strings), as they would if manipulated by a computer program. We assume that the basic operations in the structure are computable, namely, have algorithms. For a group, composition and inverse operations must be computable on these indices; for a ring or field, addition and multiplication must be similarly computable. Algorithms for other tasks on our algebraic structure can then use these basic operations as subroutines. Conversely, undecidability results in this framework can show that some natural tasks admit no uniform algorithm, and hence fail for some concrete computable systems.


The value of such a theory for developing new algorithmic methods in algebra is clear: an algorithm designed in this setting applies whenever the underlying algebraic structure has the required computable presentation. Both papers develop effective procedures for natural algebraic constructions, including substructures and extensions, and for homomorphisms and isomorphisms between represented structures. They also exhibit failures of effectiveness. Such negative results are equally important: they identify extra hypotheses needed for an algorithm and may guide the search for those hypotheses. Far less obvious is the fact, demonstrated here and in later work, that introducing computation into algebra can also yield new purely algebraic consequences that make no reference to algorithms in their final statement.

If we abstract the procedures that carry out the operations as oracles, these papers anticipate later black-box algebraic models, which likewise abstract their operations. Theories of black-box groups originated with Babai and Szemer\'edi~\cite{BabaiSzemeredi1984}, and those of black-box fields with Boneh and Lipton~\cite{BonehLipton1996}. Both focus on \emph{efficient algorithms} rather than finite algorithms, but a prescient remark at the end of Rabin's paper already suggests replacing the class of computable functions in his theory by any class of algorithms with appropriate closure properties. 


Let us conclude with several contributions distinctive to Rabin's paper~\cite{Rabin1960ComputableAlgebra}, drawn from its introduction and its two main parts.

The vision with which Rabin opens the paper is remarkable for such a young researcher shortly after completing his Ph.D. He puts the nascent theory of computability on a par with topology in its possible impact on algebra. Just as one asks whether isomorphisms, homomorphisms, and other constructions involving topological groups are continuous, Rabin proposes asking whether the analogous maps and constructions in effectively represented groups are computable. The same program is equally relevant to rings, fields, and other algebraic systems.


Part I of Rabin's paper is devoted to a general theory of computable algebraic structures, demonstrated mainly for groups. It contains results about computable homomorphisms, subgroups, and factor groups of computable groups, together with examples showing that the necessary effectiveness hypotheses cannot simply be omitted. Especially striking are the results concerning the word problem. Here the group is already given as an abstract finitely generated group, rather than solely through an arbitrary finite presentation: for a fixed finite generating set, the word problem asks whether a word represents the identity. In a computable presentation this is decidable, since one can evaluate the word. Rabin proves the converse as well.

\begin{theorem}[Rabin~\cite{Rabin1960ComputableAlgebra}]\label{thm:word-computable}
A finitely generated group has a decidable word problem if and only if it is computable.
\end{theorem}


At this point one might wonder whether noncomputable finitely generated groups exist at all. Rabin obtains them in an elegant way from the Novikov--Boone theorem and Theorem~\ref{thm:word-computable}.

\begin{theorem}[Rabin~\cite{Rabin1960ComputableAlgebra}]
There exists a finitely generated noncomputable group.
\end{theorem}

Theorem~\ref{thm:word-computable} also has a striking purely algebraic consequence.

\begin{corollary}[Rabin~\cite{Rabin1960ComputableAlgebra}]\label{cor:matrix-representation}
A finitely generated group with an undecidable word problem has no faithful finite-dimensional representation\footnote{namely, an isomorphism to a matrix group.} over any field.
\end{corollary}

Assume for contradiction that the group had such a representation, namely were isomorphic to a finitely generated subgroup of $\mathrm{GL}_n(F)$, for some finite $n$ and some field $F$. Rabin separately proves that such a group {\em is} always computable, {\em regardless of whether $F$ is computable or not}\footnote{This is a nontrivial existential argument. The key to it is that, since only finitely many elements of $F$ appear as entries of the generators, there exists a computable subfield of $F$ in which all computations can be performed.}. So, by Theorem~\ref{thm:word-computable} it has a decidable word problem, a contradiction. We note that purely algebraic examples of finitely generated groups with no faithful finite-dimensional matrix representation had previously been constructed by Fuchs--Rabinowitsch~\cite{FuchsRabinowitsch1940a,FuchsRabinowitsch1940b}; Rabin's computability-theoretic route is surprisingly direct.



Part II of Rabin's paper is devoted to computable fields. Its central positive result concerns the basic problem of constructing the algebraic closure of a field.

\begin{theorem}[Rabin~\cite{Rabin1960ComputableAlgebra}]\label{thm:computable-algebraic-closure}
Every computable field $F$ has a computable algebraic closure $\overline{F}$, and a computable embedding of $F$ into $\overline{F}$.
\end{theorem}

Rabin points out that one cannot prove this general theorem by simply carrying out the usual Steinitz construction, successively adjoining roots of polynomials through a sequence of finite field extensions. Maintaining an effective presentation along that route would require polynomial factoring algorithms over the intermediate fields, and such algorithms need not exist. Instead, Rabin constructs the closure through a different route due to Bourbaki\footnote{For experts, Rabin is working with recursively recognizable maximal ideals in a polynomial ring in countably many variables, using effective root-counting and elimination arguments.}. This route avoids the need for effective factoring algorithms.

The comparison with Fr\"ohlich and Shepherdson is especially illuminating here. They use the Steinitz construction mentioned above, and obtain a different result, in which a closure need not be computable in their sense. For their framework and definitions, a computable algebraic closure exists if and only if the required factoring algorithms exist (which sometimes do not). Rabin's theorem above proves that in his framework, a computable (in his sense) algebraic closure always exists. For background and subsequent developments, see~\cite{StoltenbergHansenTucker1999,MetakidesNerode1979}.




To sum up, Rabin's paper creates a language of computable presentations of algebraic systems and relations between them. It proves when standard algebraic constructions preserve computability, supplies counterexamples showing the limits of those preservation principles, and converts algorithmic questions, such as the word and splitting problems, into structural properties of groups, fields, embeddings, and substructures. Computability (and later, computational complexity) thus became an integral part of algebra. 

\section{Finite-state machines and the birth of nondeterminism}
\label{sec:finite-automata}

The previous section dealt with Turing machines, the most general model of computation. In this section we focus on the simplest model, \emph{finite-state machines}\footnote{We will return to Turing machines, and the boundary between decidability and undecidability, in the next section, where we will need the notions of finite automata and nondeterminism discussed in this one.}. Indeed, the most basic finite automaton, a {\em one-way deterministic finite automaton}, abbreviated $1\mathrm{DFA}$, can be viewed simply as a Turing machine devoid of two central features: the ability to write on its tape, and the ability to reverse direction. Expressed positively, a $1\mathrm{DFA}$ has a finite number of states (some of which are marked ``accepting''). It scans its input from left to right, and at every step its current state and the scanned input symbol determine the next state. Thus, the current state---a finite amount of information---is its only memory of what the automaton has seen so far. An input is accepted if the final state is an accepting state. The \emph{language} accepted by an automaton is the set of inputs it accepts.

Rabin and Scott's paper \emph{Finite Automata and Their Decision Problems} \cite{RabinScott1959} is the centerpiece of this section. It was written in the summer of 1957 at IBM Research. Both were PhD students of Church at Princeton; Rabin had just graduated, and Scott would graduate the following year. The paper earned them the 1976 A. M. Turing Award. The citation for the award summarizes its impact concisely~\cite{ACMTuringAward1976}:
\begin{quote}
\emph{This paper introduced the idea of nondeterministic machines, which has proved to be an enormously valuable concept. Their (Scott \& Rabin) classic paper has been a continuous source of inspiration for subsequent work in this field.}
\end{quote}
\noindent We shall elaborate on it.


$1\mathrm{DFA}$s were already studied extensively in the 1950s, and were extremely well understood. Kleene~\cite{Kleene1951,Kleene1956} (an early version appeared in 1951) completely characterized their computational power. He introduced \emph{regular languages}\footnote{Regular languages later featured as the simplest class of formal languages in the Chomsky hierarchy that we will mention in the next section.}, and proved that they are precisely those accepted by $1\mathrm{DFA}$s. Kleene also defined closure properties of regular languages, which led to their concise definition using \emph{regular expressions}. Myhill and Nerode gave structural characterizations of regular languages \cite{Myhill1957,Nerode1958}, and of the minimal number of states in an automaton recognizing each one. Moore and Mealy~\cite{Moore1956,Mealy1955} further studied $1\mathrm{DFA}$s that also produce an output symbol at every step.


With this state of knowledge about the limited power of finite memory, it may be hard to imagine what else can be learned about such machines. The Rabin--Scott paper demonstrated what it really means to understand a computational resource. They initiated a complete study of how different modes of computation and variants of the original model can use finite memory, and the relationships among the language-recognition power and memory usage of these different modes and variants. They also addressed the decidability of basic semantic properties of languages computed by each of these variants, as was done for Turing machines. This comprehensive treatment was foundational for the developing methodology of computational complexity theory.


\subsection{Nondeterministic computation}

By far, the most important mode of computation introduced in the paper is \emph{nondeterminism}. The state-transition rules of a nondeterministic automaton allow it, given any state and scanned input symbol, to move into any one of several possible states (instead of precisely one in a deterministic automaton). An input is accepted if \emph{there exists} a sequence of such choices that leads to a final accepting state. The resulting model is called a $1\mathrm{NFA}$, with ``N'' standing for nondeterminism, here and henceforth. Informally, a nondeterministic automaton is allowed to \emph{guess} the best possible move to make at every step. The main theorem in their paper is that, despite this powerful guessing ability, it does not enable the recognition of more languages. Abusing notation a bit and letting model names represent the classes of languages they accept, we have:

\begin{theorem}[Rabin--Scott~\cite{RabinScott1959}]
$1\mathrm{NFA}=1\mathrm{DFA}=\mathrm{REG}$.
\end{theorem}

The proof is the famous ``subset construction.'' The simulating deterministic automaton has a state for every subset of states in the simulated nondeterministic automaton. This construction entails an exponential blow-up in the number of states.


From the outset, Rabin and Scott make several key observations about this new model of computation. First, it can save memory: they state that nondeterministic automata can use far fewer states to compute the same language. They do not prove any lower bound, but a tight exponential lower bound, due to Lupanov~\cite{Lupanov1963}, soon followed. Furthermore, Rabin and Scott demonstrate how this new mode of computation (together with their theorem) is a powerful descriptive tool allowing simpler proofs and new theorems. One such simple proof establishes Kleene's closure theorems for regular languages. A new theorem we will soon discuss is the simulation of two-way automata by one-way automata. In both cases, the ability to guess is extremely powerful.

Nondeterministic computation immediately became central to formal language theory. It proved essential for exact machine characterizations of language classes in the Chomsky hierarchy, which Chomsky~\cite{Chomsky1956,Chomsky1959} developed in the 1950s as part of his theory for understanding language. In 1962, Chomsky~\cite{Chomsky1962Pushdown}
characterized \emph{context-free} languages (the second level of Chomsky's hierarchy) as those recognized by \emph{nondeterministic} pushdown automata\footnote{Here it is known that the deterministic model accepts strictly fewer languages.}. In 1964, Kuroda~\cite{Kuroda1964} characterized \emph{context-sensitive} languages (the third level of Chomsky's hierarchy) as those computed by linear-space \emph{nondeterministic} Turing machines\footnote{Here it is an open question whether the deterministic model accepts strictly fewer languages or not.}.


In parallel, nondeterminism immediately became central to logic and verification, where finite automata on \emph{infinite} objects are studied. We will discuss these at length in Section~\ref{sec:computability}, and mention here only that the central models of this kind---Büchi's automaton on infinite words and Rabin's automaton on infinite trees---are nondeterministic.

But the impact on computational complexity is astounding. Complexity theorists have lived and breathed nondeterminism ever since that paper. The natural view of nondeterministic complexity as the complexity of \emph{verification} highlights the fundamental relation between proofs and computations and its many manifestations. Many central complexity classes have an ``N'' version, and we study its effect on a variety of models and resources; we list just a few examples.

For space complexity we have Savitch's theorem~\cite{Savitch1970}, surprisingly showing that nondeterminism is not so powerful for saving memory, incurring only a quadratic saving over deterministic computation (as long as space is not too tiny, namely $s(n) \geq \log n$).
\[
\mathsf{NSPACE}(s(n))\subseteq \mathsf{DSPACE}(s(n)^2).
\]
For time complexity, however, we know almost nothing better than the trivial simulation, which incurs an exponential cost when simulating nondeterministic computation deterministically.

\[
\mathsf{NTIME}(t(n))\subseteq \mathsf{DTIME}(2^{t(n)}).
\]
For the choice of polynomial time $t(n)=n^{O(1)}$, improving this trivial simulation gives rise to the \Pclass{} versus \NPclass{} problem---among the most important open questions in computer science and mathematics. In the mid 1960s Edmonds~\cite{Edmonds1965MatroidPartition,Edmonds1967Branchings} explicitly asked\footnote{Using his terminology of ``good algorithms'' and ``good characterizations''.} about the relation between the classes \Pclass{}, \NPclass{}, \coNPclass{} and \NPclass{} $\cap$ \coNPclass{}. In the early 1970s came the foundational papers of Cook, Karp, and Levin~\cite{Cook1971,Karp1972,Levin1973} initiating the theory of \NPclass{}-completeness. These exciting developments surely had an impact on the 1976 Turing award for Rabin and Scott, and the citation for it quoted above.



The guessing power of nondeterministic programs is naturally encoded by an existential quantifier, which leads to connections with logic and suggests further classes obtained by adding a universal quantifier and allowing the quantifiers to alternate, giving rise to the polynomial hierarchy of Meyer and Stockmeyer~\cite{MeyerStockmeyer1972,Stockmeyer1976}, another fundamental family of complexity classes. Nondeterminism further leads to an alternating-quantifier characterization of \PSPACE{} in~\cite{ChandraKozenStockmeyer1981}. This relation was essential, for example, in Shamir's proof~\cite{Shamir1992} that $\mathsf{IP}=\mathsf{PSPACE}$. Indeed, the class \IP{} of interactive proofs~\cite{Babai1985,GoldwasserMicaliRackoff1989} combines nondeterminism and randomness\footnote{Such alternation is already present in Papadimitriou's ``Games against nature'' paper~\cite{Papadimitriou1985Games} which gives another characterization of \PSPACE{}.}, a development which revolutionized the notion of proofs, creating a rich theory with vast theoretical and practical applications. In yet another direction, alternation of existential and universal quantifiers was added to finite automata by Chandra, Kozen and Stockmeyer~\cite{ChandraKozenStockmeyer1981}. Drusinsky and Harel~\cite{DrusinskyHarel1994}, showed that such automata can be {\em doubly} exponentially more succinct than deterministic ones (so, another exponential is obtained by adding the universal quantifier to the existential one). More on applying such alternation (and more) to other models can be found in~\cite{HirstHarel1994,HarelRosnerVardi1990}.




And the list goes on and on.

\noindent
{\bf Remark} (Turing's choice machines)
One paragraph in Turing's famous 1936 paper~\cite{Turing1936} briefly defines what he called \emph{choice machines}, or c-machines, which allow several next moves at any given computational step. He motivates them by proof verification, and distinguishes them from his deterministic \emph{automatic machines}, or a-machines, which are the main focus of his paper. In a later footnote he shows that a finite number of choices can be eliminated in finite time by exhaustive enumeration on a deterministic a-machine\footnote{Needless to say, Turing cared only about finite algorithms and not their complexity, so the cost of this exhaustive enumeration was immaterial to him.}. This allows him to use c-machines as an intermediate step when simulating the lambda calculus, and then to convert them to a-machines by such exhaustive enumeration.

In short, it is clear that Turing (prescient as ever) defined nondeterministic computation, and viewed it as useful, long before the Rabin--Scott paper. However, this aspect of his paper was apparently not referenced by anyone for decades to come\footnote{A 1985 book review by Israel~\cite{Israel1985} identifies Turing's c-machines with nondeterministic Turing machines.}. One can speculate on the reasons for this. First, Turing's paper is so full of fundamental ideas that were crucial for the development of computer science and technology. Second, the cursory manner above in which c-machines occur in the paper does not call attention to a clear separate notion of importance. Finally, it seems that Turing never returned to develop it.


In contrast, Rabin and Scott's paper highlights this notion, provides the first systematic study of nondeterminism as a computational mode, and demonstrates its utility and power for finite space complexity. This helped make Rabin--Scott the standard source for the language theory and the logic/verification communities (from which the complexity theory community mainly grew). This may explain why the 1976 Turing Award citation credits Rabin and Scott with the introduction of nondeterministic computation.

Interestingly, almost all standard textbooks in complexity theory define the central concept of nondeterministic computation without any reference to its origins. Perhaps this demonstrates how obvious the notion has become to everyone. Given the above, I feel that Turing should be credited with the first definition of this concept, and Rabin--Scott (who, like everyone at the time, did not seem to know about it) should be credited with its reinvention, development, and impact.

\subsection{Other modes and variants of finite-state machines}

Rabin and Scott did not stop with nondeterminism. They used finite automata as a setting in which to compare several modes of computation.

\paragraph{Two-way automata.}
A $1\mathrm{DFA}$ reads the input from left to right. A $2\mathrm{DFA}$ may move its input head both left and right, like a Turing machine. This seems like a significant increase in power: the machine can revisit earlier scanned symbols, and thus use the input itself as a kind of read-only workspace. Nevertheless, Rabin and Scott prove that two-way motion does not increase the class of recognized languages.

\begin{theorem}[Rabin--Scott~\cite{RabinScott1959}]
$2\mathrm{DFA}=1\mathrm{DFA}=\mathrm{REG}$.
\end{theorem}

As noted above, the simulation proof of $2\mathrm{DFA}$s by $1\mathrm{DFA}$s is in two steps: first, a simulation of a given $2\mathrm{DFA}$ by a $1\mathrm{NFA}$, fully utilizing the power of nondeterminism, followed by a $1\mathrm{DFA}$ simulation of that $1\mathrm{NFA}$ using the theorem above. Again, the simulation can be costly in states, namely in memory. Thus one-way and two-way automata have the same decision power, but may differ sharply in memory complexity. Indeed, later work of Kapoutsis~\cite{Kapoutsis2005} proved that in this simulation there are tight exponential upper and lower bounds on the blow-up.


\paragraph{Two-way nondeterministic automata.}
Nondeterminism and two-way motion can be combined, resulting in $2\mathrm{NFA}$s. Again, Rabin and Scott prove that the resulting model still recognizes only regular languages.

\begin{theorem}[Rabin--Scott~\cite{RabinScott1959}]
$2\mathrm{NFA}=2\mathrm{DFA}=\mathrm{REG}$.
\end{theorem}

Again, the simulation of $2\mathrm{NFA}$s by $2\mathrm{DFA}$s can result in an increase in the number of states. But here the quantitative bounds are far from matching. We have an exponential upper bound, but only a quadratic lower bound on the blow-up in the number of states (see the works~\cite{SakodaSipser1978,Chrobak1986,Chrobak2003Erratum,AdeogunKapoutsis2026}). Proving a superpolynomial lower bound remains a challenging problem\footnote{This would not contradict Savitch's theorem above, as here space $s(n) = O(1)$ whereas there it must be $s(n) \geq \log n$.}.


\paragraph{Multitape automata.}
In all models considered so far, finite automata had only one tape. Rabin and Scott now introduce and study finite automata with two or more tapes. Here the story changes dramatically. The addition of multiple tapes can create genuine increases in power. In particular, two-way two-tape automata cannot in general be simulated by one-way two-tape automata. Indeed, the paper gives an explicit non-regular language\footnote{This is really a \emph{relation} on input pairs on the two tapes that are accepted. A modern discussion of this issue can be found in~\cite{Furia2012}.} which \emph{can} be accepted by a two-way automaton with two tapes! This provides an important counterpoint to the equivalence theorems: some changes in the finite-state model preserve decision power, while others do not.


\paragraph{Probabilistic automata.}
In a later paper~\cite{Rabin1963ProbAutomata}, Rabin explores yet another novel mode of computation, adding \emph{random transitions} to finite-state automata. Being a major initial work in the study of randomness as an algorithmic tool and computational resource, we defer a fuller discussion to Section~\ref{sec:probabilistic-algorithms}.

\subsection{Decision problems about automata}

The title of the Rabin--Scott paper emphasizes not only finite automata, but their decision problems. So far we have described their work on the decision problems \emph{solved} by finite automata. Their comprehensive treatment now turns to decision problems \emph{about} finite automata: problems where the input is a finite automaton, and we seek information about the language it computes. Analogous to similar work for Turing machines, Rabin and Scott turn to consider which natural semantic properties of finite automata are decidable. Typical examples include \emph{emptiness}, which asks whether the language computed by a given automaton is empty, and \emph{equivalence}, which asks whether two given automata compute the same language.



As above, there is a sharp transition between the difficulties of such problems when they concern one-tape finite automata and those concerning multitape ones. On the one hand, the paper proves that for all models of one-tape finite automata mentioned above, these and related problems are decidable. Indeed, this tractability makes these models amenable to analysis, allowing a useful theory of these models. Of course, these results build on the extensive understanding of the structure of one-tape finite automata developed before and in their paper, and the equivalences they proved between them. On the other hand, they show that even the emptiness problem for two-tape, two-way $\mathrm{DFA}$s is undecidable.

\begin{theorem}[Rabin--Scott~\cite{RabinScott1959}]
It is undecidable to decide whether a given deterministic two-tape, two-way finite automaton accepts any input.
\end{theorem}


This part of the paper connects back to the previous section and, as we shall soon see, forward to the next one. Rabin and Scott continue to explore the fine line between decidability and undecidability, and its dependence on delicate differences between the exact details of their input structures.

For further reading on automata and formal-language theory, see Eilenberg~\cite{Eilenberg1974}, Hopcroft and Ullman~\cite{HopcroftUllman1979}, Kozen~\cite{Kozen1997}, and Sakarovitch~\cite{Sakarovitch2009}.
\section{Between decidability and undecidability: logic, trees and games}
\label{sec:computability}


This section is devoted mainly to Rabin's paper~\cite{Rabin1969}, which proves the decidability of an important logical theory. In several late interviews, Rabin describes this paper as the most difficult work he ever did \cite{ShashaRabin2010,RabinHarel2015}. We will now slowly develop the background needed to explain Rabin's result, and then discuss some of its consequences. Different perspectives on this evolution can be found in the surveys of Thomas \cite{Thomas1997,Thomas1995}, Grädel--Thomas--Wilke \cite{GradelThomasWilke2002}, and Vardi--Wolper \cite{VardiWolper1986}.

\subsection{Background}

Recall that Hilbert's \emph{Entscheidungsproblem} calls for an algorithm to determine whether a given mathematical statement is valid \cite{HilbertAckermann1928}. Of course, to make the question formal, we must fix an underlying proof system.


The most basic logical theory is the so-called Peano arithmetic. It is a first-order theory of the natural numbers with addition and multiplication, with familiar axioms and deduction rules. Many theorems in large parts of mathematics are proved in this system. It was thus especially devastating that, in the negative answer to Hilbert's \emph{Entscheidungsproblem}, the work of Gödel, Church, and Turing proved that, for this very theory, testing validity is undecidable~\cite{Godel1931,Church1936Entscheidungsproblem,Turing1936}. The same holds for testing the {\em provability} of a given statement\footnote{By Gödel's \emph{completeness} theorem, validity and provability are the same for most proof systems in use, but anyway their undecidability proof applies to both without using that theorem. When talking about the undecidability of a proof system, we would usually mean that testing validity is undecidable, but may switch to testing provability.}. 



Indeed, Turing proved undecidability for any \emph{sufficiently strong} logical theory, where this term informally means that the theory can argue about computation. For any such theory, there is no algorithm to decide whether a given formula is valid, or if it is provable. This naturally leads to the question: \emph{Which logical theories are decidable?} We will see in this section how different changes in either the logic or the underlying structure can completely change the answer.


In 1929, Presburger~\cite{Presburger1929} proved that the first-order theory of the natural numbers with addition but \emph{without} multiplication is decidable. This so-called Presburger arithmetic is quite strong and can prove many statements about ``linear structures'' within the integers. But the absence of multiplication allowed Presburger to develop an algorithm for the validity of any formula. Note that this result predates Gödel's incompleteness theorem, and was a positive result within Hilbert's program.

Another major decidability result was obtained in 1951 by Tarski~\cite{Tarski1951}. He proved that the first-order theory of the real numbers with addition and multiplication is decidable. Thus, replacing the integers by the reals restores decidability even while retaining both field operations. To prove his result, Tarski developed ``quantifier elimination,'' a central algorithmic technique.

The study of the two theories above was motivated by mathematical considerations: understanding Diophantine problems about natural numbers and about systems of polynomial equations and inequalities in real-closed fields, respectively. The study of the following theories was motivated by computer-science considerations. This program, connecting logic and language theory with computation, was initiated by Church~\cite{Church1957CircuitSynthesis} with his \emph{analysis} and \emph{synthesis} questions. The analysis question asks us to determine the properties, expressed logically or syntactically, of the input--output behavior of computer programs. The synthesis question asks us to design a program whose input--output behavior meets specifications so expressed. For both, we seek the relation between the expressive power of the logical theory describing the properties and the computational power of the programs or machines that realize them. This study has remained important throughout the history of computer science, as computational models evolve and ensuring that what is specified is what is implemented---in hardware, software, protocols, and now AI models---remains essential to the field.


This program was met with spectacular success in the following years within formal language theory (which we mentioned in the previous paragraph and is worth repeating here). The first major success of this program was due to Kleene~\cite{Kleene1956} (who also coined the terms ``analysis'' and ``synthesis''). He defined \emph{regular languages} and proved that they are precisely the languages computed by one-way deterministic finite automata. Then Chomsky~\cite{Chomsky1962Pushdown} characterized the richer class of \emph{context-free languages}\footnote{Essential for the developing theory of programming languages.} as those computed by \emph{nondeterministic pushdown automata}. And then Kuroda~\cite{Kuroda1964} characterized the even richer class of \emph{context-sensitive languages} as those computed by nondeterministic linear space Turing machines. In all cases, the equivalences (analysis and synthesis) were algorithmically constructive.



But another main focus in this area, in work following Church in the late 1950s by Elgot, Trakhtenbrot, and Büchi~\cite{Trakhtenbrot1958,Elgot1961,Buchi1960}, was the development of different frameworks for specifying desired properties of computations, and studying the computational models they capture. One difference is that the specifications in the above results are {\em language-theoretic} (a grammar for the language), while these new works used {\em logical} specifications (which we will describe after discussing the models). Another difference is the emphasis on properties that hold throughout the computation. Many computer systems\footnote{And, with hindsight, also biological systems (for example, animal behavior) and societal systems (for example, emergency services).} interact with their environment: operating systems, for example, repeatedly receive input signals from processes that request a resource (or raise an alarm or effect a change), and respond with some output action. For such systems, especially ones we design, we typically impose requirements on the relation between input and output that should be satisfied throughout.


A natural abstraction of this very general scenario is called, in modern language, a \emph{reactive system},  defined in this generality by Harel and Pnueli~\cite{HarelPnueli1985}\footnote{In the original papers such systems have different names, including sequential circuits, logical nets, and Moore automata. A very related modern term (typically with different applications in mind) is the notion of an \emph{online algorithm}.}. In reactive systems, the inputs $I$ are words over some alphabet, and the responses $R$ are words over the set of states of our system. At time step $t$, an input symbol $I_t$ arrives, and the system, given its current state, processes it and moves to a new state, $R_t$, in a manner consistent with its transition rules. We may now specify, for example, using logical formulas, the sets of input--output pairs of words $(I,R)$ that satisfy the requirements. A possible output $R$ is usually called a \emph{run}, as it describes the behavior of the system as we ``run'' it on input $I$.



What should be the weakest formal theory needed to specify natural properties of pairs $(I,R)$? The discussion below hints at the art of navigating the delicate balance between the expressiveness of a logical theory and the complexity of deciding it. Clearly, a first-order formula over the integers can specify consistency of the local transition rules at every step (for both deterministic and nondeterministic computation). For finite computations, it can check that the final state is accepting. Moreover, as the only ``arithmetic'' in action is repeatedly shifting time ahead by one, we do not even need addition (as in Presburger arithmetic), and can instead allow only the ``add one'' successor function, $s(t)=t+1$.

However, first-order statements in this setting do not suffice to check most interesting properties of computations. Basic properties of reactive systems as above, for example, ``every request is eventually satisfied,'' require quantification over \emph{subsets} of time steps. First-order theories allow only quantification over elements of the underlying structure (which for time steps, can be taken to be the integers). Theories allowing quantification over predicates over the universe are called \emph{second-order}, and are vastly more powerful. Still, an important observation was that most desired properties\footnote{For example, that no error ever occurs, that no processor ``starves'' forever, that every packet is eventually delivered, etc.} need only quantifiers over \emph{unary} predicates, namely over subsets of the universe. Such theories are called \emph{monadic second-order} theories, which are still very strong. Can we expect decidability of such strong theories?

%
The decisive positive answer was given in Büchi's 1962 paper~\cite{Buchi1962}. He proved decidability of \SoneS, the monadic second-order theory over the integers with the successor function. 

\begin{theorem}[Büchi~\cite{Buchi1962}]
\SoneS\ is decidable.
\end{theorem}

Büchi's proof has several remarkable features. First, he notes that, as \SoneS{} can argue about subsets of integers, it can argue about infinite words. This leads him to extend the usual model of finite automata (whose inputs are finite words) to operate on infinite words. Specifically, {\em Büchi automata} are one-way, \emph{nondeterministic}\footnote{This choice is essential, and in the setting of infinite words deterministic Büchi automata are strictly weaker. However, McNaughton~\cite{McNaughton1966} later showed that nondeterministic Büchi automata can be simulated by deterministic \emph{Muller} automata, a variant of Büchi's. This turned out to be important in Rabin's proof. An exponentially faster determinization was achieved by Safra~\cite{Safra1988}; this paper had major impact on the field of verification.} finite automata whose inputs are \emph{infinite} words. When does such an automaton accept an input? For a finite input, we accept if the final state is accepting. But for an infinite input, the automaton works forever and never reaches a final state. A Büchi automaton is defined to accept an input $I$ (an infinite word) if some possible run $R$ on $I$ enters an accepting state infinitely often. Observe that this condition can easily be expressed in \SoneS. Büchi automata are sometimes called $\omega$-automata, and languages over infinite words are sometimes called $\omega$-languages.



Next, Büchi proves an equivalence between logic and automata. To every \SoneS-formula there corresponds a Büchi automaton such that the formula is valid if and only if the automaton accepts \emph{all} inputs. More generally, the satisfying assignments of the formula can be made to be exactly the language accepted by the automaton. The proof of this theorem proceeds by a natural induction on the structure of the formula. The second-order existential quantifier corresponds to projection, negation is complementation, and the connectives, conjunction and disjunction of formulas become intersections and unions of languages. By far, complementation is the hardest to achieve, as can be expected from the nature of the accepting condition\footnote{Interestingly, this argument uses Ramsey theory in an ingenious way which allows the complement of a language to become a \emph{finite} union of (simpler) languages.}.



Finally, and this is the main novelty in Büchi's paper, the equivalence above is used to prove his computability theorem. The algorithm for testing the validity of an \SoneS-formula first constructs the associated automaton, and then checks that it accepts all inputs. This idea was a methodological transformation affecting logic, automata theory, and program verification. We note that the complexity of this algorithm, while finite, is nonelementary, which Meyer~\cite{Meyer1975} proved is unavoidable.



\subsection{Branching computation and Rabin's paper}
This result of Büchi became one of the foundations of the now vast fields of verification and model checking, capturing many important properties of \emph{linear} computation, like safety, liveness and fairness. But even Büchi himself was not completely satisfied, and asked what could be done about \emph{branching} computation \cite{Buchi1962}. In such a computation, the input resides on the vertices of, say, a binary tree, rather than on a linear path. A possible run is again a labeling of the vertices by states of the system. Computation by an automaton is still local: transition rules relate the input $I_v$ and the state $R_v$ at any node $v \in \{0,1\}^*$ to the states labeling its two children, $R_{v0}$ and $R_{v1}$. The accepting condition should be satisfied for every infinite path down from the root.


Observe that a theory describing properties of such branching computations would be far more expressive. Specifically, it would be able to address, at once, all possible behaviors of the environment interacting with our system (whereas in linear systems the environmental behavior is encoded in the input). In other words, it allows specification over all possible futures. The analogous logical theory, the monadic second-order theory of \emph{trees} (as opposed to words), is a simple modification of \SoneS: one allows \emph{two} successor functions instead of one, $s_0(v)=v0$ and $s_1(v)=v1$, and it is appropriately called \StwoS. The decidability of this powerful theory was posed as an open question by Büchi in his paper \cite{Buchi1962}.

Rabin solved this problem.
\begin{theorem}[Rabin~\cite{Rabin1969}]
\StwoS\ is decidable.
\end{theorem}

As with Büchi, Rabin proves this via an equivalence between formulas and automata. To every formula in \StwoS, Rabin constructs a finite automaton over binary trees, which yields a correspondence between satisfying assignments of the formula and inputs accepted by the automaton. This \emph{Rabin tree automaton} is again a nondeterministic finite automaton, and on any input labeling tree vertices, a possible run labeling them must satisfy local consistency of the transition rules for every vertex and its two children. What becomes far subtler is defining an accepting condition. Rabin puts a condition on every branch (infinite path down from the root of the tree). Fixing finitely many pairs of subsets of the automaton's states, $(E_i,F_i)$, \emph{Rabin's acceptance condition} demands that along every branch there is an index~$i$ such that the run visits some state in~$F_i$ infinitely often, but visits every state in~$E_i$ only finitely many times. This delicate choice is critical for the equivalence and the algorithm. Both follow the structure of Büchi's proof, namely induction on formula structure, but implementing them is far more complex. Indeed, to accept an input, there must exist a run in which every branch is accepting. This requires a very different complementation procedure for which Rabin invents new tools. His paper is considered a technical tour de force.



The impact of Rabin's paper was enormous, and we conclude by discussing some of it. The insightful definitions and new techniques Rabin developed became essential tools in the field of verification and model checking. Most consequences of the main result are due to the strength of a theory on branching structures. Rabin's paper already demonstrates this by giving, as corollaries, new decision procedures for an impressive list of logical theories (including ordered, arithmetic, topological, and lattice-theoretic theories), by embedding them all in \StwoS. The expressivity of the branching structure greatly expands the set of computational properties that can be captured by a decidable theory. As mentioned, branching allows us to consider at once all eventualities caused by the environment (including, for example, adversarial behavior, timing issues, orderings of actions by multiple processes, and more), and formulas can express relationships between different branching futures. Rabin's theorem influenced many developments of other logics, such as modal and temporal logics, and decision procedures for them. Surveys illuminating various aspects of this impact include \cite{Thomas1997,Thomas1995,GradelThomasWilke2002,VardiWolper1986}.



\subsection{Infinite games and determinacy}

We conclude with another advantage of having tree structure: it allows us to model and study infinite games. At a conceptual level, the verification and model checking literature often views synthesis as a game played between the program and the environment, where a winning strategy for the program means meeting the given specifications. But Rabin also shows concrete corollaries for the theory of infinite games. This theory, initially developed by Gale--Stewart~\cite{GaleStewart1953} and continued by Wolfe~\cite{Wolfe1955} and Martin~\cite{Martin1975}, extends \emph{determinacy} theorems to larger and larger classes of infinite games. These results extend Zermelo's determinacy theorem for finite win--lose two-player games~\cite{Zermelo1913}, which states that some player always has a winning strategy. Rabin's paper shows how Wolfe's theorem follows easily as a corollary from the decidability of \StwoS. Interestingly, Rabin had already made the following important point about this theory in an earlier paper~\cite{Rabin1957}. He gives simple examples of infinite games that are determined, namely, one of the parties always has a winning strategy, but none of that player's winning strategies is computable. Thus, the utility of having a winning strategy is doubtful, and this result calls for investigating the existence of computable ones when defining determinacy.




\section{Computational complexity}
\label{sec:complexity}

Computational complexity is today one of the most exciting and influential intellectual disciplines, with great impact not only on computer science, but also on mathematics, the exact and social sciences, and technology. Part of its success depended on extremely prescient choices made when it was born in the 1960s, choices that were central in allowing its initial structures, goals, techniques, and results to emerge. These were then refined and expanded repeatedly to yield the theory and applications we have today.

The centerpiece of this section is Rabin's 1960 paper, \emph{Degree of Difficulty of Computing a Function and a Partial Ordering of Recursive Sets}~\cite{Rabin1960Difficulty}, in which he initiated a general computational complexity theory. In it he lays out and motivates, in clear language, many of the important choices and focal points that we now take for granted. We discuss these foundational issues first. When reading this part, try to imagine a world when quantitative analysis of algorithms hardly existed.

\begin{remark}
We note that this paper was only published as a technical report of the Hebrew University. However, we know that these ideas have influenced later papers that were foundational to complexity theory, including those of Cobham~\cite{Cobham1965}, Ritchie~\cite{Ritchie1963}, and Hartmanis--Stearns~\cite{HartmanisStearns1965}. Cobham cites Rabin's paper in his. Ritchie's paper acknowledges conversations with Rabin. Hartmanis, in his memoir on the origins of computational complexity theory~\cite{Hartmanis1981}, tells how talking to Rabin when he and Stearns were writing their paper sent them to ``eagerly read his manuscript.''

\end{remark}


\subsection{Basic notions, choices and goals}

\medskip
\noindent\textbf{Complexity matters.}
Church, Turing, and Post initiated modern computability theory by establishing a fundamental fact about computation: some natural problems are uncomputable. But this is only a starting point. Not all computable functions are born equal, and having a finite algorithm may be nearly as useless in practice as having none. Rabin opens his paper with an example of two functions, one that is trivial to compute (it takes linear time), and one that seems much harder (the only known algorithm takes exponential time and is unrealistic for inputs of reasonable size). He initiates the research program of formalizing the notion of the ``degree of difficulty of computing functions,'' which we now call computational complexity.

\medskip
\noindent\textbf{General complexity measures.}
Rabin points out that time may not be the only resource we want to study; we may wish to consider other complexity measures that signify difficulty of computing functions. He thus proceeds to define complexity measures axiomatically. 

\medskip
\noindent\textbf{Axiomatic, machine-independent approach.}
An important consideration when choosing axioms is machine independence. Even if we focus on a fixed resource, say time, the measured complexity of a problem depends on various details of the computational device, which would obscure rather than illuminate. Rabin therefore seeks conclusions that do not rest on one particular implementation. Rabin proposes a framework to achieve this, defining complexity measures using three axioms, and his hierarchy results (below) then hold for \emph{every} measure satisfying these axioms. 

Blum~\cite{Blum1967} later built on Rabin's work, and developed a somewhat different, simpler, and more streamlined axiomatic framework for the same purpose, which became standard. Blum complexity measures are defined via two simple axioms, which are completely machine independent. Besides a hierarchy theorem for his measures, Blum additionally proved fundamental compression and speed-up results for them.


\medskip
\noindent\textbf{Decision problems.}
To avoid extraneous effects on complexity, such as the size of an output, Rabin focuses his theory on functions taking only the values $0$ and $1$, namely decision problems.

\medskip
\noindent\textbf{Asymptotic complexity.}
Rabin suggests that we should ignore the complexity of any finite set of inputs, and rather focus on eventual, limiting behavior of any algorithm. Fix a class of algorithms and a complexity measure $m$. We say that $f<g$ (read ``$g$ is harder than $f$'') if there is an algorithm $A$ computing $f$ such that, for every algorithm $B$ computing $g$, one has
\[
 m_A(n)<m_B(n)
\]
for all but finitely many $n$. 


\medskip
\noindent\textbf{The partial order of complexity.}
For any admissible complexity measure, Rabin studies the partial order it induces on computable functions as defined above. He proves that this partial order contains infinite chains, and in particular functions of arbitrarily high difficulty for any such measure.

\begin{theorem}[Rabin~\cite{Rabin1960Difficulty}]\label{thm:rabin-complexity-hierarchy}
For every computable function $f$, there is a computable function $g$ taking only the values $0$ and $1$ such that $f<g$. If $f$ is primitive recursive, $g$ may also be chosen primitive recursive.
\end{theorem}

For concrete resources, stronger and more detailed class-theoretic results became possible. Hartmanis and Stearns~\cite{HartmanisStearns1965} developed explicit time-bounded classes, simulation results, and hierarchy theorems. A closely related follow-up by Hartmanis, Lewis, and Stearns~\cite{HartmanisLewisStearns1965} classified computations and proved hierarchy results for both time and space requirements. By now we have numerous such results.




We now turn back to focus on time complexity, and on the major choice of what the new theory of computational complexity should treat as an \emph{efficient} or \emph{tractable} algorithm.

\subsection{The emergence of the class \texorpdfstring{\Pclass{}}{P}}

The choice of polynomial time algorithms as the theoretical counterpart of practically efficient ones was crucial for the field. Four influential early papers advocating this choice are those of Edmonds~\cite{Edmonds1965}, Cobham~\cite{Cobham1965}, and Rabin~\cite{PerlesRabinShamir1963,Rabin1967Automata}. Their respective arguments and emphases are overlapping, and together make a strong case for this choice for the notion of {\em tractability}. We summarize their arguments below. 

Cobham stresses the robustness of polynomially bounded computation across reasonable machine models (for example, various types of Turing machines). Cobham and Rabin emphasize that the closure of polynomials under addition and multiplication means that polynomial bounds on complexity are stable under common programming practices of composing computations. Edmonds and Rabin draw motivation from concrete algorithms: Edmonds from his algorithm for perfect matching in general graphs, and Rabin from his algorithms with Perles and Shamir~\cite{PerlesRabinShamir1963} for the synthesis and analysis of automata. In both cases these algorithms have polynomial running times of low degree, while straightforward exhaustive procedures were exponential. In these early days before complexity considerations became widespread, Edmonds and Rabin felt compelled to explain why the former are much better than the latter, as exponential growth is much faster. Rabin adds a vivid way of illustrating it: if computers become faster by a constant factor, a fixed amount of time permits a multiplicative increase in input size for a polynomial-time algorithm, but only an additive increase for an exponential-time algorithm.


These arguments proved convincing, and the community adopted polynomial time as the central initial formalization of efficiency and tractability. This choice led directly to some of the most important complexity classes and questions. The class \Pclass{} contains problems with efficient deterministic algorithms; \BPP{} contains problems with efficient probabilistic algorithms and bounded error; \NPclass{} contains problems with efficiently verifiable certificates. Moreover, the \NPclass{}-complete problems become the hardest problems in \NPclass{}, under efficient reductions. These are only the most prominent examples: polynomial resource bounds enter the definitions of hundreds of classes and their relationships; see, for example, the Complexity Zoo~\cite{ComplexityZoo}. 

The impact is broader still. Modern theories of cryptography and pseudorandomness replaced particular attacks and statistical tests with the class of {\em all} polynomial-time adversaries or distinguishers. The compositional closure properties of this class again allowed researchers to uncover structure and greatly develop these theories. These are specific examples of a more general principle, suggesting that basic beliefs pervading philosophy about knowledge had to radically change. The old (information theoretic) idea that you know everything which you can infer from your current knowledge was usefully and practically replaced in many settings by the new (complexity theoretic) idea that this holds only for what you can {\em efficiently} derive (e.g., in \Pclass{}) from your current knowledge (see the survey of Aaronson~\cite{Aaronson2013} for an excellent survey).


One cannot exaggerate the importance of this initial choice. In essence, the flexibility of polynomial time gave complexity theory \emph{robust} notions with which to explore relationships among problems and models. These revealed many structures, such as complexity classes, as well as flexible techniques that created a working methodology of the field. As this understanding evolved, the theory could go on to study many more refined computational models and resources, leading to more refined results and questions, and to closer ties with practical efficiency and relevant technologies (as well as to beautiful mathematics).

%

\subsection{Lower bounds}

From its inception, computational complexity aimed to prove that some problems are difficult to compute. Rabin~\cite{Rabin1960Difficulty} opens his paper by asking whether inverting von Neumann's function (which for an $n$-bit input $x$ outputs the $n$ ``middle'' bits of $x^2$) is hard. Cobham~\cite{Cobham1965} starts (and ends) his paper by asking whether multiplication is harder than addition. Of course, we do not seem any closer to resolving the \Pclass{} versus \NPclass{} problem than we were when it was formulated in the early 1970s.


Seven decades of attempts to prove hardness have taught us that this is among the hardest tasks in computational complexity. Progress along this front has been slow and has focused on restricted models of computation, with the hope of developing stronger and stronger techniques that may allow progress. Rabin has been part of this effort, proving mainly lower bounds on arithmetic and algebraic circuit models, for example, in his papers~\cite{Rabin1972Positivity,Rabin1972ScalarProducts,RabinWinograd1972,BaurRabin1982}, which we do not review here. We conclude by describing a different important early contribution of Rabin to proving lower bounds. 

In 1974, Fischer and Rabin~\cite{FischerRabin1974} proved a powerful lower bound on the complexity of any decision procedure for Presburger arithmetic. Recall from Section~\ref{sec:computability} that this is the first-order logical theory of the integers, which allows addition but not multiplication. It was proved decidable by Presburger in 1929, but the complexity of his algorithm was nonelementary\footnote{Namely, requiring time growing faster than any fixed tower of exponentials.}. How hard is it? Fischer and Rabin show that it requires at least doubly exponential time\footnote{Single-exponential time lower bounds had been given for some natural problems in the 1972 paper of Meyer and Stockmeyer~\cite{MeyerStockmeyer1972}, and also in the present paper of Fischer and Rabin for deciding the first-order theory of the reals.}.

\begin{theorem}[Fischer--Rabin~\cite{FischerRabin1974}]
Any algorithm deciding the validity of formulas in Presburger arithmetic must run for at least $2^{2^n}$ time on formulas of size $n$.
\end{theorem}


Their proof is a very sophisticated diagonal argument. To implement it, they encode the computation performed by any algorithm $A$ on input $x$ of length $n$ as a formula in Presburger arithmetic, which is valid if and only if the algorithm halts in $2^{2^n}$ time. In essence, the formula of length $n$ has to somehow encode the tableau of the entire computation, and check its correctness. Doing so for such a long computation---for example, accessing exponentially long addresses and extracting the bit stored at a given address---seems impossible in a very succinct formula with the severely restricted expressive power of this logic. Thus the heart of their argument is precisely inventing such an encoding.

To date, no better lower bound has been proved. Oppen~\cite{Oppen1978} greatly improved the upper bound by giving a decision procedure that runs in triply exponential time. The single exponential gap remains an open problem.
\section{Probabilistic algorithms}
\label{sec:probabilistic-algorithms}

%
%

The field most associated with Rabin is that of probabilistic algorithms. Today this is a huge field, with numerous such algorithms, analysis tools, and application areas, which are the subject of many textbooks and standard courses. Rabin played a major and distinctive role in its rapid early development in the 1970s. A sequence of his papers (and lectures) modeled, exemplified, and enthusiastically advocated the use of randomness (and the toleration of error) in computation. They demonstrate how probabilistic algorithms can achieve better efficiency in various resources, as well as greater simplicity in their design and analysis, for many problems in diverse domains. In all of these, Rabin insisted from the start on the now-standard worst-case analysis over inputs, of the resources used and the error incurred by these algorithms.

This section describes some of his probabilistic models and algorithms, which in turn improve memory use, running time, communication, and simplicity over deterministic ones. In subsequent sections on cryptography and distributed computing, we will see more of his probabilistic algorithms and protocols, which outperform deterministic ones in different ways. We will also discuss Rabin's work on modeling probabilistic computation and reasoning about it.


\subsection{Background}

Randomness has both mystified and been utilized by humans since ancient times. Many societies have used lotteries of different sorts for decision making and for gambling. Mathematical studies of randomness started with probability theory (Pascal, Laplace), and proceeded through statistics (Galton, Pearson, Fisher), statistical mechanics (Boltzmann, Maxwell), algorithmic randomness (Kolmogorov, Solomonoff, Chaitin), information theory, and private-key cryptography (Shannon). The use of randomness as a novel mathematical tool arose in the probabilistic method (von Neumann, Shannon, Erd\H{o}s).

The computing revolution brought the use of randomness into computation. Turing was perhaps the first to explicitly suggest using a ``random element'' in (learning) algorithms, in his famous paper on artificial intelligence~\cite{TuringAI}; he even suggested physical devices that would generate such randomness. Early uses of randomness in computation were primarily in simulations of natural phenomena, such as weather prediction and nuclear interactions, and in many other applications mentioned in the famous Metropolis--Ulam paper~\cite{MetropolisUlam1949} expounding the \emph{Monte Carlo} method. It was also used for numerical integration. These initial works were mostly for practical uses, without a formal model or performance analysis. More fundamentally, many of these applications were statistical in nature, where the need for randomness is almost explicit in the problem formulation, and for which deterministic algorithms seemed irrelevant. In contrast, most problems we shall consider in this section are discrete problems from different domains, such as decision or search problems. Classically, these were solved only by deterministic algorithms, and allowing randomness (which may have seemed irrelevant) offered surprising gains.



\subsection{Modeling probabilistic algorithms}

The first paper to suggest a formal probabilistic model of computation is by de Leeuw, Moore, Shannon, and Shapiro~\cite{DeLeeuwMooreShannonShapiro1956}, who defined a probabilistic Turing machine in 1956. However, the authors study its power from a computability (rather than computational complexity) viewpoint. Their main result is that this model is not more powerful than the usual deterministic Turing machines; i.e., such machines can compute with high probability only decidable functions. The simulation is by the usual enumeration of all possible values of the random choices, and counting how many of them lead to acceptance. Thus, for decidability, there is no advantage in using probabilistic algorithms over deterministic ones.


The focus on computational complexity changes everything. The first paper to do so was written a few years later, in 1963, by Rabin~\cite{Rabin1963ProbAutomata}. It defines a different formal model, namely probabilistic finite automata. In the spirit of his paper with Scott discussed in Section~\ref{sec:finite-automata}, it continues to investigate finite space complexity. He proves, as with nondeterminism, that probabilistic algorithms can provably save space. As usual, Rabin's paper provides a comprehensive treatment of the new model, its variants, properties, and their relationships. Here we describe the main import of this paper. In the next subsection we discuss his 1976 paper modeling and exemplifying time-efficient probabilistic algorithms.

Probabilistic finite automata are defined as a natural extension of deterministic and nondeterministic finite automata: for every current state and symbol read, the transition function gives a probability distribution on the next state. The machine accepts an input if the probability of reaching an accepting state exceeds a threshold $\lambda$.

From the outset, Rabin distinguishes two types of machines:
\begin{itemize}
\item ``unbounded error,'' in which an input is rejected when its acceptance probability is at most $\lambda$; and
\item ``bounded error,'' (or ``gapped'') in which, for some \emph{gap} $\delta>0$, an input is rejected when its acceptance probability is at most $\lambda-\delta$, with a guarantee that no input has acceptance probability in $[\lambda-\delta,\lambda]$.
\end{itemize}
Analogously to the notation of Section~\ref{sec:finite-automata}, we denote a gapped one-way probabilistic finite automaton by $1\mathrm{PFA}$\footnote{Rabin calls the threshold $\lambda$ a \emph{cut-point}, and the bounded-error model one with an \emph{isolated cut-point}.}.


Following Rabin, this distinction was adopted, became central, and has since been investigated in complexity theory for a variety of computational resources; we mention some. For example, for probabilistic polynomial-time computations, we have the class \PP{} of unbounded-error algorithms and the class \BPP{} of bounded-error algorithms. In this model, we have several reasons to believe that the unbounded-error model is more powerful\footnote{On the one hand, \PP{} contains \NPclass, which is believed to require exponential-size circuits. On the other hand, \BPP is known to have polynomial-size circuits. A separation of \PP{} and \BPP also follows if the polynomial-time hierarchy does not collapse, combining~\cite{Toda1991} and~\cite{Sipser1983,Lautemann1983}.}. For other models and resources, like space complexity, communication complexity, and decision-tree depth, we can formally prove (without any assumptions) that the unbounded-error model is more powerful.

%

For finite automata, Rabin proved that this distinction is essential. On the one hand, he proves that the unbounded error model accepts undecidable languages (giving both a counting argument proof as well as an explicit one)\footnote{This result necessitates arbitrary real numbers as threshold. If one demands rational (or even a computable real) threshold, and rational transition probabilities, the resulting language is decidable.}. On the other hand, he proves that the bounded error model computes only regular languages. Ever conscious of computational resources, Rabin also gives upper and lower bounds on the blow-up in the number of states (the memory requirement) in simulating probabilistic machines by deterministic ones.


%

\begin{theorem}[Rabin~\cite{Rabin1963ProbAutomata}]
\hfill
\begin{itemize}
\item $1\mathrm{PFA}=1\mathrm{DFA}=\mathrm{REG}$.
\item An $n$-state, $\delta$-gapped $1\mathrm{PFA}$ can be simulated by a $(1/\delta)^{O(n)}$-state $1\mathrm{DFA}$.
\item For every $\delta$ there exists a $\delta$-gapped two-state $1\mathrm{PFA}$ such that any $1\mathrm{DFA}$ simulating it requires $1/\delta$ states.
\end{itemize}
\end{theorem}

In the proof of the second item, the states of the $1\mathrm{DFA}$ simulating a given $1\mathrm{PFA}$ are an appropriate $\epsilon(\delta)$-net of probability distributions over the states of the $1\mathrm{PFA}$, with the natural transition function between them. This may be viewed as a quantitative generalization of the famous subset construction in the Rabin--Scott simulation of $1\mathrm{NFA}$s by $1\mathrm{DFA}$s.

Rabin proceeds to study a natural computational question suggested by this simulation, namely approximating the product of a given sequence of stochastic matrices. He points out that the naive approach cannot be space-efficient, even if the entries are small rationals, because of the numerical precision needed as the sequence grows longer. He then proposes a natural condition on the inputs under which a direct reduction to the simulation above actually works in finite memory. This basic problem was later studied in detail, as it was naturally recognized as a complete problem for simulating probabilistic log-space machines; see, for example,~\cite{BorodinCookPippenger1983,SaksZhou1999,Hoza2022} for the history and best-known results. The paper further studies how other aspects of a model as simple as finite automata are affected by the injection of randomness, including the stability of probabilistic automata under perturbations of transition probabilities, as well as special types of such automata; we shall not elaborate on these.

Instead, we conclude with later work on variants of this model, which further illustrates its subtleties. Curiously, Rabin does not study $2\mathrm{PFA}$s, namely bounded-error probabilistic automata that, as in deterministic automata, allow two-way reading access to the input. Their power was discovered by Freivalds~\cite{Freivalds1981TwoWay}, who showed that they recognize the nonregular (but context-free) language $\{0^n1^n:n\ge 1\}$. This is surprising as it informally means that, with high probability, {\em fixed} constant space suffices to count to $n$ for any {\em unbounded} integer $n$.

This clever algorithm (which the reader is challenged to discover) runs in exponential time in the input length. Perhaps even more surprisingly, Dwork and Stockmeyer~\cite{DworkStockmeyer1990} proved that, when restricted to polynomial expected time, $2\mathrm{PFA}$s can recognize only regular languages. Models that combine randomness and nondeterminism in finite automata were studied as well, and their exact power is still a mystery (see~\cite{CondonHellersteinPottleWigderson1998} for partial progress and an intriguing conjecture).


\subsection{Primality testing: randomness for algebra and number theory}

In April 1976 Joe Traub organized a meeting at Carnegie Mellon University titled ``Algorithms and Complexity: New Directions and Recent Results.'' The invited speakers included many leaders of the already quite diverse field at the time (including Berlekamp, Brent, Garey, Johnson, Karp, Knuth, Kung, Paterson, Rabin, Shamos, Simon, Tarjan, Traub, and Yao). The proceedings of the meeting~\cite{Traub1976} present quite a panorama justifying the title: the volume displayed the remarkable breadth of the exciting new trends and directions then emerging. Rabin was there to excite the audience about probabilistic algorithms. Excitement was needed: Rabin felt compelled to state explicitly in his paper, ``It may at first seem surprising that employing randomization leads to efficient algorithms,'' a statement no one would dream of uttering a decade later.
Rabin's paper, ``Probabilistic Algorithms''~\cite{Rabin1976ProbabilisticAlgorithms}, is the focus of this subsection, and the next one. But first, some background.


The 1970s saw probabilistic algorithms become a central algorithmic tool to be considered for any problem that deterministic algorithms were used to solve, in the hope of beating them! Some famous examples existed even beforehand. A unique occurrence, long before computers were dreamed of, is Pocklington's 1917 paper~\cite{Pocklington1917} on computing square and cube roots in finite prime fields. While lacking precise language and a formal model, the paper clearly gives a polynomial-time probabilistic algorithm for these problems (which to this very day do not have a deterministic counterpart without number theoretic assumptions!). This paper remained generally unknown, but the idea was rediscovered several times (for example, by Peralta, Lehmer, and Berlekamp~\cite{Peralta1986,Lehmer1969,Berlekamp1970}, as we shall see below). Another famous example is Hoare's 1961 paper~\cite{Hoare1961,Hoare1962} introducing his probabilistic sorting algorithm \emph{Quicksort}. While it has no asymptotic advantage over the $O(n \log n)$ deterministic sorting algorithms such as \emph{mergesort}, Hoare argues that his algorithm is simpler, easier to implement, and faster by a constant factor in expectation.


An important paper in the development of probabilistic algorithms is Berlekamp's 1970 paper~\cite{Berlekamp1970} on factoring univariate polynomials over large prime fields. Polynomial factoring is a central, basic algebraic problem. It generalizes root taking (above) and plays an important role in decoding algebraic codes such as Reed--Solomon codes. Importantly, we do not know any deterministic polynomial-time algorithm for this problem.

\begin{theorem}[Berlekamp~\cite{Berlekamp1970}]
There is a probabilistic expected-polynomial-time algorithm that, given a polynomial $f(x)$ and a prime $p$, produces the (unique) factors of $f$ in $\mathbb{F}_p[x]$.
\end{theorem}


But perhaps the largest impact demonstrating the algorithmic power of randomness came from the probabilistic \emph{primality-testing} algorithms of Solovay--Strassen~\cite{SolovayStrassen1977} and Rabin~\cite{Rabin1976ProbabilisticAlgorithms,Rabin1980Primality} in the late 1970s.


\begin{theorem}[Solovay--Strassen and Rabin~\cite{SolovayStrassen1977,Rabin1976ProbabilisticAlgorithms,Rabin1980Primality}]
There is a polynomial-time probabilistic algorithm for testing whether a given integer is prime or composite. It is always correct on prime inputs and errs with low probability on composite inputs.
\end{theorem}


It may be important to point out a source of power of randomness in such algebraic and number-theoretic settings, which may differ from that in other applications. Without describing any of the algorithms in the theorems above, let me note how algebraic ``rigidity'' easily supplies the error bounds for them. Essentially, random choices are bad (lead to error) only if they fall into a ``small'' algebraic substructure of the structure they are chosen from (for example, a proper subgroup or coset, or an algebraic subvariety of a vector space). In the polynomial-factoring algorithms and the primality testing algorithms, the random choices are typically elements of a group $G$ (for example, $\mathbb{Z}_N^{\ast}$, the multiplicative group of invertible elements modulo $N$ for some integer $N$). Simple analysis shows that a random choice leads to an error if and only if it lies in a proper subgroup $H$ of $G$ (e.g., containing elements of a trivial multiplicative character value). Then, Lagrange's elementary theorem that the size of a subgroup always divides the size of the group directly shows that the error has probability at most 1/2! \footnote{Another important set of examples of how algebraic structure is utilized is when the random choices are taken from a vector space, and lead to error  when a low-degree polynomial vanishes on them. For univariate polynomials (a vector space of dimension 1, namely one variable), the degree of the polynomial trivially upper-bounds the number of roots (on which it vanishes). This leads to a simple error bound when the number of possible choices is larger than the degree. For multivariate polynomials the same error bound holds as well, due to the important but simple-to-prove multivariate extension in~\cite{DeMilloLipton1978,Schwartz1980,Zippel1979} of the univariate case. This principle underlies many important probabilistic algorithms, including the celebrated algorithm for \emph{multivariate} polynomial factoring over finite fields due to Kaltofen~\cite{Kaltofen1986,Kaltofen1989} (here both input and outputs are given by straight line programs).}




The impact of these two (independent but very similar) primality tests was huge. They announced the arrival of probabilistic algorithms as an independent field, causing everyone seeking algorithms to consider employing randomness just as naturally as seeking a deterministic algorithm. Part of the appeal and inspiration of the probabilistic primality tests must have come from the fame of the problem itself (for example, even Gauss appealed to the mathematics community to find efficient ways to test primality\footnote{And to factor integers.}). Another attraction must have come from how essential large primes became to public-key cryptography with the foundational works of Diffie--Hellman, Rivest--Shamir--Adleman~\cite{DiffieHellman1976,RSA1978}, Goldwasser--Micali~\cite{GoldwasserMicali1982,GoldwasserMicali1984}, and the wave of theory and practice that followed.

Indeed, the period of the late 1970s and the 1980s saw an explosion of probabilistic algorithms in numerous domains. These included polynomial identity testing~\cite{DeMilloLipton1978,Schwartz1980,Zippel1979}, universal hashing~\cite{CarterWegman1979}, randomized routing~\cite{Valiant1982Routing}, parallel matching algorithms~\cite{KarpUpfalWigderson86,MulmuleyVaziraniVazirani1987}, randomized rounding~\cite{RaghavanThompson1987}, and many others.


Let us return to Rabin's 1976 paper~\cite{Rabin1976ProbabilisticAlgorithms}, as it contains far more than the primality-testing algorithm. Rabin clearly set out to excite that conference audience, by laying a foundation for a new and exciting subfield of algorithms, and arguing for its importance and potential with new examples and open questions. Again, today many of these are considered standard and obvious, but apparently this was not so at the time.

Rabin begins his paper by distinguishing these algorithms from another subfield, probabilistic analysis of algorithms (the subject of Karp's lecture at the same conference), in which one analyzes \emph{average-case} behavior of deterministic algorithms under some natural input distributions. Rabin stresses that probabilistic algorithms demand \emph{worst-case} bounds on their complexity and error bounds. He lays out the model in some detail. A complete formalization and the familiar classes followed a year later in Gill's paper~\cite{Gill1977}. He discusses the variants that we today call $\BPP, \mathsf{RP}, \mathsf{ZPP}$: respectively, two-sided-error, one-sided-error, and zero-error models. Only then does Rabin demonstrate the power of randomization by giving two novel efficient probabilistic algorithms from very different domains: number theory (the primality test) and geometry (which we discuss in the next section). These two algorithms use very different techniques.


We proceed with the primality story we started earlier. Rabin notes that his algorithm is really a version of a deterministic primality test of Miller~\cite{Miller1976},\footnote{Indeed, it is often referred to as the Miller--Rabin test.} which runs in polynomial time under the ERH (the Extended Riemann Hypothesis, a generally believed conjecture in number theory; it guarantees that good random choices exist in a small interval, which may be searched deterministically). Rabin suggested that ERH may be false, and that a random search will find good random choices with high probability. Another point Rabin makes is that his (and Solovay--Strassen's) one-sided algorithm certifies composites with high probability but does not certify primes. Given that short certificates for primes do exist (as Pratt showed in~\cite{Pratt1975}), Rabin asks for a one-sided algorithm in the other direction, which would find such certificates with high probability. Such algorithms were eventually found, first by Goldwasser--Kilian~\cite{GoldwasserKilian1986} (under mild number-theoretic assumptions) and then by Adleman--Huang~\cite{AdlemanHuang1992} (without them). These algorithms and their analyses are much subtler and surprisingly use the theory of elliptic (and higher genus) curves. Combined with the original algorithms, one gets zero-error primality tests running in expected polynomial time, which always certify their answer. The beautiful ending of this story is, of course, the deterministic polynomial-time algorithm for primality by Agrawal, Kayal, and Saxena~\cite{AgrawalKayalSaxena2004}. It is worth mentioning that it emerged from ingeniously \emph{derandomizing} yet another probabilistic primality test of Agrawal and Biswas~\cite{AgrawalBiswas2003}, which is very different from those above. The best version of this deterministic algorithm runs in about $O(n^5)$ time, and to this day most practical systems (for example, in cryptography) use the much faster, roughly $O(n^2)$-time original probabilistic primality tests! 



Rabin continued looking for algebraic and number-theoretic problems to exemplify the power of randomness. In joint work with Shallit~\cite{RabinShallit1986}, he gives a probabilistic algorithm for finding a representation of a given integer as a sum of four squares; the existence of such representations is guaranteed by a famous theorem of Lagrange~\cite{Lagrange1770}. In~\cite{Rabin1980FiniteFields}, Rabin gives probabilistic algorithms for several problems in finite fields, including finding irreducible polynomials, finding roots of univariate polynomials, and factoring them.



\subsection{Closest pairs: randomness for geometry}

Another subfield of algorithms that burst into existence in the 1970s was computational geometry. Among several initial works, a major impetus came from Shamos's Ph.D. thesis~\cite{Shamos1978Thesis}, papers, and his survey of 75 problems~\cite{Shamos1974Problems}. Basic problems in which the typical input was a collection of points in Euclidean space were addressed, many arising from practical questions about the real world. These problems differ from the general combinatorial and graph-theoretic questions on which the algorithmic community was focused, as geometry imposes extra structure, including convexity, order, and metric structure. Thus, completely new questions could be asked, for example, about convex hulls and Voronoi diagrams, while old questions such as shortest paths and minimum spanning trees in weighted graphs were restricted to the special case in which weights came from Euclidean distances between points. New techniques were developed to design and analyze algorithms for such problems. 


To demonstrate the power of randomness in this setting, Rabin picks the \emph{closest-pair problem}: given a set of $n$ points in the Euclidean plane, find the pair that is closest. The best deterministic algorithm for this problem, due to Hoey and Shamos~\cite{ShamosHoey1975}, runs in $O(n \log n)$ time\footnote{This algorithm actually solves a more general problem, of finding a closest point for every point in the set.}. Moreover, like many geometric algorithms, it relies on sorting, and so this running time seems the best one can hope for. Rabin shows that with randomness it can be solved in linear, $O(n)$ time.


\begin{theorem}[Rabin~\cite{Rabin1976ProbabilisticAlgorithms}]
There is an expected $O(n)$-time probabilistic algorithm that finds the closest pair among $n$ points in the plane. The same holds in $\mathbb{R}^d$ for any fixed dimension $d$.
\end{theorem}



Unlike primality testing, where both the algorithm and its analysis are simple and require a few lines each to describe, here both are considerably more complex and subtle\footnote{Even when I oversimplify and ignore the precise machine model.}. With this example Rabin shows that the use of randomness can be far richer than one-shot sampling that is guaranteed to land in a large set. Let me briefly sketch the main ideas behind his solution.

The algorithm is recursive and demonstrates the power of subsampling. Given a set $P=P_0$ of $n$ points, we select a subset $P_1$ of $n/2$ of them at random, then a subset $P_2$ of these of size $n/4$, etc. The idea is to use the solution to $P_i$ to solve the problem on $P_{i-1}$ in linear time. This leads to a time recurrence of the form $t(n) = t(n/2) + O(n)$, yielding $t(n) = O(n)$.

To facilitate the inductive step, Rabin proves that a strong invariant is maintained with high probability. If the shortest distance between a pair of points in $P_i$ is, say $2\delta_i$, then a grid of this resolution (specifically of side length $\delta_i$) will have at most one point of $P_i$ (and, in particular, the closest pair will be in neighboring squares in this grid). Now add the points of $P_{i-1} \setminus P_i$ back to this grid. Then solve the closest-pair problem in the neighborhood of each nonempty cell by brute force, comparing all pairs in it. This is a quadratic amount of work per cell, and so the crux of the analysis is proving that the randomness of these sets guarantees, with high probability, that the sum of these squared occupancies is still $O(n)$.



Needless to say, a full description of the algorithm and the necessary data structures to support it efficiently is considerably more complex (indeed, some of these details are missing, and were supplied later: see the paper of Banyassady and Mulzer~\cite{BanyassadyMulzer2018} for full details and simpler analysis).


Subsampling and recursive subsampling became basic tools in probabilistic algorithms. One beautiful example is the linear-time probabilistic algorithm for minimum spanning trees of Karger, Klein, and Tarjan~\cite{KargerKleinTarjan1995}, which has a very similar high-level structure.


\subsection{Fingerprinting: randomness and communication}

Another powerful demonstration of randomness gives provably exponential savings in communication. The basic setup, laid out in Yao's seminal paper~\cite{Yao1979} on communication complexity, is the following. Two players, Alice and Bob, hold $n$-bit strings $x$ and $y$, respectively, and must compute a fixed function $f(x,y)$ using a fixed protocol of their choosing. The communication complexity of $f$ is the minimum number of bits they must exchange to determine the function value for every input.

Easily the most basic function to consider is equality: testing whether $x=y$. In lectures and classes, Rabin motivated this problem, as well as the focus on communication costs and ignoring computational ones, through what he called the ``man on the moon'' problem. Alice on Earth and Bob on the moon want to verify consistency of the two huge data files they hold. Of course, one of them can send their whole input to the other, at a possibly prohibitive cost of $n$ communicated bits. This unfortunately is best possible; it is not hard to see that in any deterministic protocol they must exchange at least $n$ bits. The Rabin--Yao probabilistic protocol\footnote{It is referred to by Yao~\cite{Yao1979} as joint unpublished work with Rabin in the paper above, and is attributed to both in all standard texts~\cite{MotwaniRaghavan1995,KushilevitzNisan1997,MitzenmacherUpfal2005}.} does exponentially better!


\begin{theorem}[Yao~\cite{Yao1979}]
There is an $O(\log n)$-communication probabilistic protocol for equality with inverse-polynomial error.
\end{theorem}


The protocol and its analysis are extremely simple. Alice picks a small random prime $p \leq n^2$, and sends Bob both $p$ and $x' = x \pmod p$, only $O(\log n)$ bits in total. Bob computes $y' = y \pmod p$ and decides according to whether $x'=y'$. Clearly, if $x=y$, this decision is always correct. If $x \neq y$, at most $n$ distinct primes\footnote{Indeed, at most $n/\log n$ primes, if one uses standard results on the distribution of primes.} can divide $x-y$; thus, the error probability is at most $O(1/n)$.

The very same idea was independently discovered and used at about the same time for a related, though different, demonstration of the power of randomization for one-tape Turing machines. It was known~\cite{Barzdin1965,Hennie1965} that such machines require $\Omega(n^2)$ time to recognize the equality language, namely input strings of the form $xx$. Freivalds~\cite{Freivalds1977} uses the same reduction modulo a small random prime to show that randomness gives a near-quadratic saving.

\begin{theorem}[Freivalds~\cite{Freivalds1977}]

There is a probabilistic one-tape Turing machine that recognizes the equality language with high probability in time $O(n\log n)$.
\end{theorem}


This idea of \emph{random fingerprinting} is fundamental and in universal use! Indeed, it is so basic that it is regularly taught in undergraduate classes without giving a reference for its origin. Needless to say, people used reduction modulo primes as simple {\em hash}\footnote{A synonym for fingerprint.} functions in practice, but typically used fixed primes. Freivalds, Rabin, and Yao all stress that random fingerprinting works for \emph{every} input with high probability; despite compressing large objects into tiny fingerprints, collisions rarely occur!

Rabin kept seeking different methods and different applications of random fingerprinting. In his paper~\cite{Rabin1981Fingerprinting}, Rabin suggested that, instead of viewing bit strings as encoding integers and reducing them modulo primes (which requires modular arithmetic), one can view them as $\mathbb{F}_2$-polynomials and reduce them modulo irreducible low-degree polynomials (which requires only simple local bit operations). A different method, encoding bit strings by random $\mathrm{SL}_2(p)$ matrices, is proposed and analyzed in his paper with Karp~\cite{KarpRabin1987}, which is the subject of the next subsection.


\subsection{Pattern matching: randomness and simplicity}

Rabin's next demonstration of the algorithmic power of randomness concerns a central problem that already had an optimal algorithm. Here the power is to \emph{simplify}: to replace complex algorithms by a simpler one with a simple analysis, broad generalizability, and better practical efficiency.

\emph{Pattern matching} is the following basic problem in text processing. Given two sequences (say, binary, although any alphabet will work), a pattern $x$ (of length $m$) and a text $y$ (of length $n \geq m$), the goal is to find occurrences of $x$ as a consecutive block in $y$.

A trivial $O(mn)$ algorithm suggests itself. However, Knuth--Morris--Pratt~\cite{KnuthMorrisPratt1977} and Boyer--Moore~\cite{BoyerMoore1977} found ingenious $O(n)$-time algorithms\footnote{All algorithms in this subsection are in the RAM model, which allows constant-time operations on logarithmic-size words.} for this problem. The Karp--Rabin paper \emph{Efficient Randomized Pattern-Matching Algorithms}~\cite{KarpRabin1987} applied the random fingerprinting method above to give a probabilistic algorithm that runs in expected $O(n)$ time. Despite not improving the time complexity, they make the important point that simplicity (in exposition, analysis and/or implementation) may sometimes justify the use of randomness and the tolerance of error. Moreover, their algorithm is far more space-efficient than the deterministic one.


\begin{theorem}[Karp--Rabin~\cite{KarpRabin1987}]
There is an $O(n)$-time, $O(\log n)$-space probabilistic algorithm for pattern matching with negligible error.
\end{theorem}

Again, the idea is strikingly simple. As pattern matching calls for a sequence of equality tests (of $x$ to length-$m$ ``windows'' of $y$), fingerprinting is natural to use. Thus, we fingerprint $x$ by a logarithmic-size value $x'$ as above, similarly fingerprint the first $m$-window of $y$ by $y'_1$, and compare the fingerprints. The main observation is that, as the window slides one symbol to the right, updating the fingerprint $y'_i$ to $y'_{i+1}$ requires only a \emph{constant} number of additions, subtractions, and comparisons, resulting in a linear-time algorithm\footnote{If we want to find, for example, the first occurrence of $x$, we can even make the algorithm zero-error and still run in expected linear time by checking the putative match (investing time $m$ whenever $x'=y'_i$), as the probability of a false match can be reduced below $1/(mn)$ by picking a random prime below, say, $(mn)^2$.}.


Note that in this algorithm (and indeed, the deterministic ones), the text $y$ can be given as a stream rather than being presented fully before the computation begins. Unlike the deterministic algorithms, which require memory $\Omega(m)$ bits, the probabilistic one uses only $O(\log n)$ memory. Streaming algorithms with low memory requirements are by now a huge field (see, for example,~\cite{Muthukrishnan2005}), and the Karp--Rabin algorithm can be seen as a very early example. For more details on this viewpoint see~\cite{BreslauerGalil2014}.

\subsection{Logic for probabilistic computation}

Rabin's interest in randomness was not limited to algorithm design. Once probabilistic computation becomes a central paradigm, one also needs formal systems for reasoning about it.

In joint work with Halpern, Rabin studied logics for reasoning about likelihood and probabilistic computation \cite{HalpernRabin1987}. This work belongs conceptually with the logical themes of the earlier sections, but it is naturally discussed here because the computations being reasoned about are probabilistic.

The point is that probabilistic computation introduces new kinds of assertions. One may wish to say not only that an event occurs, but that it occurs with high probability, or that one event is more likely than another. Such statements require logical languages richer than ordinary Boolean assertions.

For this purpose, Halpern and Rabin introduce a new qualitative modal logic $LL$ for probabilistic computation, and compare it to other proposals for that purpose. For $LL$, they gave possible-worlds semantics, proved a complete axiomatization, and gave an exponential-time decision procedure for satisfiability.


This completes a circle in Rabin's work. In the decidability section we saw logic used to reason about computation and automata. Here, after randomness has become an algorithmic tool, logic is again brought in to reason about the new computational phenomenon. Further information about subsequent works and this general field can be found in the monograph~\cite{Halpern2017Uncertainty} and the survey~\cite{DemeyKooiSack2013}.


%
%
%
%

\section{Cryptography}
\label{sec:cryptography}

Cryptography is another vast subject; again, we introduce it only briefly, and we refer the reader to the textbooks of Goldreich~\cite{Goldreich2001,Goldreich2004} and of Katz and Lindell~\cite{KatzLindell2014}, or the much shorter cryptography section in~\cite{Wigderson2019}.

Since time immemorial, people have had secrets that they wished to share with some and keep private from others. Thus, schemes for secret communication were developed in very ancient cultures. What was common to all of them, from Caesar's substitution code to the Enigma machine and everything in between, was the assumption that the parties wishing to communicate secretly already shared some secret information. Moreover, this assumption was formally justified by Shannon, who proved, using his newly developed \emph{information theory}, that the parties' shared secret key must be completely random and at least as long as the message they send if they want no information about it to be leaked to an eavesdropper~\cite{Shannon1949}\footnote{The alert reader must have noticed the circularity here: to communicate secretly at some point, one must have communicated secretly at some earlier point, and so on. Shannon proves that this cannot be avoided if perfect information-theoretic secrecy is required.}.


The revolution that completely removed the need for this ``essential'' assumption came in 1976 with the breakthrough paper of Diffie and Hellman~\cite{DiffieHellman1976}, which introduced computational complexity into cryptography. Their paper proposed making the entirely reasonable assumption that all participants are computationally limited 
(to polynomial time) and, furthermore, that some problems cannot be solved so quickly. In particular, they proposed the \emph{discrete logarithm} problem as such a hard problem. Under this assumption, they showed that secret communication is possible between any two parties sharing no prior secrets whatsoever. To understand how this circumvents Shannon's barrier, we must understand a fundamental difference between the information-theoretic and complexity-theoretic settings. In the former, every representation of a piece of information is equivalent to every other. In the latter, this symmetry can be broken, as some representations are easier to compute than others. 

Fix a prime~$p$ and a generator~$g$ of the multiplicative group~$\mathbb{Z}_p^\ast$. Then consider the (bijective) exponentiation map
\[
 x\longmapsto y=g^x\pmod p.
\]
Obtaining~$y$ from~$x$ is easy, using repeated squaring of $g$ modulo~$p$, requiring only polynomial time. However, obtaining~$x$ from~$y$ is the discrete logarithm problem, which is assumed intractable. Thus, in a computationally limited world, these two representations are very different. Diffie and Hellman exploit this asymmetry to create \emph{public-key cryptography}, in which any two parties sharing no prior secret information can efficiently exchange secret messages that cannot be decoded by any other computationally bounded party. They also give a completely new cryptographic application: personal digital signatures that others cannot forge. 

The main impact of the Diffie--Hellman paper was suggesting a new paradigm for cryptographic security: to prove that your protocol (for any particular task) is secure, prove that breaking it will result in an efficient algorithm for a clean, hopefully well-studied mathematical problem. One such problem, perhaps {\em the} problem, was integer factorization, which has been studied for millennia (in particular, Gauss famously appealed to the mathematical community to find a fast algorithm for it). And it seems hard: the best known classical factoring algorithm is by Buhler, Lenstra and Pomerance~\cite{BuhlerLenstraPomerance1993}; it requires time about $\exp(n^{1/3})$ to factor $n$-bit integers\footnote{In 1994 Shor~\cite{Shor1994} discovered a poly-time {\em quantum} algorithm, so if quantum computers become a reality the hardness of factoring can no longer serve cryptography; indeed, practical security systems today are being replaced by other, hopefully ``quantum resilient'' assumptions.}.




The hardness of factoring soon became the assumption of choice in theory and practice. The first to introduce it were two papers, by Rivest-Shamir-Adleman and by Rabin. Rivest, Shamir, and Adleman~\cite{RSA1978} introduced their famous RSA public-key cryptosystem, whose security is based on a problem related to the complexity of factoring integers, while Rabin introduced a public-key cryptosystem whose security is \emph{equivalent} to the complexity of factoring~\cite{Rabin1979DigitalizedSignatures} (both systems provide digital signatures as well). These works still fell short of being able to securely encrypt single bits or to argue about partial leakage of information. The seminal 1982 paper of Goldwasser and Micali~\cite{GoldwasserMicali1982,GoldwasserMicali1984} introduced \emph{probabilistic encryption} and laid formal mathematical foundations for the field, giving precise definitions of security and a public-key system that attains them based on the hardness of another problem closely related to integer factoring.


The decade following the Diffie--Hellman paper saw a frenzy of activity in the newly formed community of modern, complexity-based cryptography, exploring which other apparently impossible tasks were possible in this new wonderland. It was a period in which one could ``ask for the moon'' and not be entirely surprised to get it. Thus were born protocols for such cryptographic primitives and goals as \emph{coin flipping into a well, poker playing over the telephone, zero-knowledge proofs}, and many other feats. Together with new techniques for analyzing protocols, and new models and assumptions, the road was paved for the Internet revolution that followed in their wake.

In the next two subsections we discuss two important contributions of Rabin to this activity. We note in passing that Rabin, whose primary passion was randomness and its utility, was naturally drawn to cryptography, a field in which randomness is essential.

\subsection{Oblivious transfer and secret exchange}

Of the many ingenious cryptographic primitives invented in this hectic period, one of the most surprising and, as we shall see, fundamental is Rabin's \emph{oblivious transfer}~\cite{Rabin1981OT}. We have two parties, conventionally called Alice and Bob. Alice holds a secret. An oblivious transfer is a protocol in which Alice and Bob participate, achieving two goals: first, Bob receives Alice's secret with probability exactly one-half; and second, Alice learns nothing about whether Bob received it or not.

My first reaction---and probably that of many others---on hearing this for the first time was, ``Why?'' Who in the world would want such a functionality? In what situations would Alice prefer \emph{not} to know whether her secret went through?

Rabin, of course, knew the answer. He first constructed a simple oblivious-transfer protocol based on the hardness of factoring. He then applied this seemingly useless primitive to design a protocol for one of the challenging tasks at the time: \emph{secret exchange}. In this task, Alice and Bob each have a secret, and we seek a protocol (with extremely high probability) ensuring that Alice gets Bob's secret \emph{if and only if} Bob gets Alice's secret\footnote{Rabin makes an assumption of a physical nature, which was later removed in a sequence of works starting with Blum~\cite{Blum1983SecretExchange} on what became known as ``fair secret exchange'' using``gradual disclosure,'' in which OT played a central role as well.}.


\begin{theorem}[Rabin~\cite{Rabin1981OT}]
Oblivious transfer implies secret exchange.
\end{theorem}


This demonstration of the power of oblivious transfer was the first in a long sequence. Soon after Rabin's paper appeared, Even, Goldreich, and Lempel defined a variant they called ``1-out-of-2 oblivious transfer'' (which is equivalent in power to Rabin's) and showed further applications, some of which we shall see below~\cite{EvenGoldreichLempel1985}. A full understanding of that power of OT came later, when Kilian~\cite{Kilian1988} proved that oblivious transfer is complete for \emph{secure multiparty computation}.

\begin{theorem}[Kilian~\cite{Kilian1988}]
Oblivious transfer implies secure multiparty computation.
\end{theorem}

Secure multiparty computation is, in a sense, the mother of all cryptographic protocols. It allows a set of mutually distrustful parties to evaluate jointly \emph{any} efficiently computable, possibly randomized function~$f$ of their secrets. Indeed, being able to do this allows one to solve a vast class of cryptographic tasks. Yao~\cite{Yao1982,Yao1986} defined this problem and solved it for the two-party case. His protocol starts with a Boolean circuit for~$f$. One party first ``garbles'' the circuit and then the second party evaluates that garbled circuit. Both are carried out in a gate-by-gate fashion, according to the order of the circuit computation, where the garbling of each gate uses the oblivious transfer primitive. Goldreich, Micali, and Wigderson extended Yao's work to secure multiparty computation by any number of players~\cite{GMW1987}. Their protocol is different; all parties jointly evaluate successive shared encodings of intermediate gate values (again using OT), which are unintelligible to any individual party. These results assume an underlying public-key encryption.




Oblivious transfer is a special case of secure multiparty computation, which seems much easier. Kilian's theorem gives a converse: oblivious transfer is sufficient for general secure multiparty computation (without assuming public-key encryption), making Rabin's apparently strange primitive one of the central building blocks of cryptographic protocol design~\footnote{We note that the relative power of the primitives oblivious transfer and public-key encryption is not known, but there is no way to get one from the other via ``black-box" reductions~\cite{GertnerKannanMalkinReingoldViswanathan2000}.}. 


\subsection{Public randomness and beacons}

Complexity-based cryptography is partly the art of inventing assumptions that allow one to implement a variety of cryptographic primitives and applications. Even at the inception of the field, a distinction was made between primitives requiring only \emph{one-way functions} and those requiring the stronger assumption of \emph{trapdoor functions}\footnote{For example, it took a long time before \emph{digital signatures}, originally based on trapdoor functions, could be implemented assuming only one-way functions~\cite{Rompel1990}. Likewise, the powerful \emph{zero-knowledge proofs} require only one-way functions~\cite{GoldreichMicaliWigderson1991}.}. Since then, numerous finer or very different assumptions have been made to satisfy newly required functionalities. Needless to say, making the weakest possible assumptions is desirable, and a whole subfield of cryptography is devoted to studying their relationships.


A striking early contribution to this art was Rabin's proposal for a very different type of assumption: having \emph{public randomness} as a cryptographic resource. In his work on transaction protection by beacons~\cite{Rabin1983Beacons}, Rabin suggested a public service that broadcasts fresh random values to all participants. This assumes trust in a central entity, but the trust is minimal: it is limited to the supply of bits, independently of their use in any protocol. The idea was that a common source of unpredictable public randomness could simplify some cryptographic protocols and enable some that are impossible otherwise. In his paper, Rabin devised protocols in this setting for three related but different problems arising in digital communication, which are easy to implement when people meet in person or use trusted parties:
\begin{itemize}[leftmargin=1.5em,label={}]
\item \emph{Contract signing.} Alice and Bob each have a copy of a contract~$C$. A protocol must ensure that Alice receives Bob's signature on~$C$ if and only if Bob receives Alice's signature on~$C$\footnote{This looks superficially similar to the \emph{secret exchange} problem above because of the simultaneity requirement, but the transferable legal evidence makes the two functionalities different.}.
\item \emph{Confidential disclosure.} Alice holds a secret~$D$, which she is willing to share with Bob on the condition that he not divulge it. A protocol must ensure that Bob receives~$D$ if and only if Alice receives a legally binding statement from Bob that he will keep it secret.
\item \emph{Certified mail.} As with paper mail, Alice wants to send a message~$M$ to Bob and have its receipt certified. A protocol must ensure that if Bob receives~$M$, then Alice obtains a certificate that Bob received it.
\end{itemize}



This idea has become increasingly important. Public randomness appears in numerous protocols for coin flipping, electronic voting, noninteractive proof systems, distributed protocols, and modern blockchain systems; see, for example, direct continuations of Rabin's beacon proposal and later public-randomness services~\cite{Muftic1990Antennas,JiwaSeberryZheng1994,LeeChangChan2005,FischerIorgaPeralta2011,SytaEtAl2017}. In different applications, the assumed properties of the available public randomness vary---even Rabin's paper uses two different forms---but the essential structure remains that all parties have access to this resource. In some applications, public randomness has been proved essential\footnote{Often under assumptions such as random oracles.}, whereas in others, it has later been shown to be removable. The latter demonstrates the power of the paradigm as much as the former: it allows one to solve a problem gradually, first using public randomness and then removing it\footnote{The same has happened many times with randomized algorithms: a probabilistic algorithm is introduced first, and later a derandomization makes it deterministic. Indeed, as we mentioned in Section~\ref{sec:probabilistic-algorithms} the famous AKS deterministic primality test~\cite{AgrawalKayalSaxena2004} is an extremely clever \emph{derandomization} of the probabilistic Agrawal--Biswas test~\cite{AgrawalBiswas2003}.}.


Rabin's three applications of public randomness above provide one example of the latter. In the aforementioned paper, Even, Goldreich, and Lempel employ their version of Rabin's oblivious transfer to remove the need for public randomness in these applications~\cite{EvenGoldreichLempel1985}.


Our final substantive subsection demonstrates a completely different application of public randomness.

\subsection{Bounded storage and everlasting security}

In an insightful paper, Maurer~\cite{Maurer1992} suggested replacing the standard assumptions in complexity-based cryptography that limit the time available to participants in cryptographic protocols with an assumption limiting their storage space. His \emph{bounded-storage} model further assumes public randomness, conceptually similar to Rabin's \emph{beacon}~\cite{Rabin1983Beacons}. Specifically, we assume a stream of random bits, available to all participants, that is much longer than the given space bound. This paper initiated a beautiful line of work that achieves information-theoretic security for some cryptographic tasks under these natural assumptions. We informally describe the evolution of ideas and results, in which Rabin played a major role.


Maurer considers secret communication. Recall that, for information-theoretic security, Shannon already proved that if Alice wants to send Bob an $n$-bit message that looks completely random to an eavesdropper Eve, then Alice and Bob must share in advance a completely random secret key of length~$n$. Maurer shows that, in his model, they can share far shorter keys and still obtain complete information-theoretic security, except on a vanishing-probability event. Thus, this model provides a way to overcome Shannon's limit. Maurer also makes the point that this security holds even if the original shared secret key is later revealed. However, Maurer makes a strong restrictive assumption about the way Eve can use her memory: she can store only a subset of the bits from the public randomness. He raises the question of removing this assumption: can Eve process the public randomness as it is streamed and store any bounded amount of information about it?


In a sequence of papers with his former and current students~\cite{AumannRabin1999,AumannDingRabin2002,DingRabin2002}, Rabin removed this assumption and did considerably more. In particular, they:
\begin{itemize}
\item proved that Eve may store arbitrary information about the public randomness, provided her storage capacity is sufficiently smaller than the length of that randomness;
\item proved that Alice and Bob can use the same key for exponentially many message exchanges, achieving the same information-theoretic security as Maurer;
\item extended the applications beyond secret communication to authentication and non-malleability; and
\item formalized \emph{everlasting security}, which Maurer noted above. They proved that their protocols guarantee that after the protocol is over, even if Eve obtains the shared secret key used by Alice and Bob and her space bound is lifted (namely, she has unbounded storage), the messages sent retain the same information-theoretic security.
\end{itemize}


The papers above did not achieve optimal parameters. Near-optimal parameters were obtained later by Dziembowski and Maurer~\cite{DziembowskiMaurer2004}. A more general understanding of this model came with Lu's paper~\cite{Lu2004}, connecting the bounded-storage model with \emph{randomness extractors}. This was followed by a paper of Vadhan~\cite{Vadhan2004}, which constructed appropriate \emph{locally computable} extractors and allowed most existing protocols and results to be generalized, simplified, and given improved parameters.

\section{Distributed computing}
\label{sec:distributed-computing}

Distributed computing captures any situation in which many computers or processes aim to achieve a variety of goals (common or separate) in a way that necessitates interaction. This, of course, encompasses numerous situations in computer science and technology (as well as in the natural sciences), making distributed computing a vast and heterogeneous field. General overviews of this field can be found, for example, in the books of Lynch~\cite{Lynch1996} and Attiya and Welch~\cite{AttiyaWelch2004}, or the much shorter section on distributed computing in \cite{Wigderson2019}. We note that distributed computing is similar to cryptography in many ways, as both involve many interacting parties. However, the issues central to each are typically quite distinct, as we shall see\footnote{Naturally, there is a nontrivial intersection and plenty of interaction between these fields: cryptography is used in many distributed algorithms, and distributed algorithmic ideas are exploited in cryptographic protocols. The new area of blockchains perfectly exemplifies this mixing.}.

One focus of many distributed systems and problems is \emph{locality}. A distributed protocol assigns to every computer in the system its own program, and there is no central entity controlling them all. Each such local program iterates steps (or rounds) of internal computation (usually viewed as free) and communication exchanges with its neighbors. The latter can take place through various media: shared variables, communication channels in a given network topology (which is often not known to the individual processes), and so forth.

Besides locality, distributed protocols have to defend against faults. These can come in a variety of forms. Obvious examples include failures of, or even adversarial behavior by, some processors or communication media. Another difficulty, which will appear in some of the examples below, arises from \emph{asynchrony}. In synchronous systems, we assume a global clock, and so the rounds of internal computation and message exchange with neighbors are synchronized. But in many realistic situations one cannot assume a global clock, and the local clocks may bear no relation to one another. This asynchrony is typically modeled by an adversary that has a global view of the system and can use it to schedule processor actions and message delays arbitrarily. Even in synchronous settings, we let such an all-knowing adversary dictate which processes fail and how, often dynamically during execution. A distributed protocol has to achieve its goals typically against such powerful adversaries (whose precise nature may differ and is specified for each problem). This adversarial power explains the difficulty of designing distributed protocols and the subtleties involved in proving their properties\footnote{For this reason it would be impossible to completely elaborate the full assumptions on the models and adversaries; these can be found in the references.}.


In this section we survey major contributions made by Rabin to four fundamental problems (of very different natures) in distributed computation: \emph{dining philosophers}, \emph{mutual exclusion}, \emph{Byzantine generals}, and \emph{information dispersal}. The issues they raise are as relevant to practical applications today as they were when formulated decades ago. One need only consider the Internet and modern distributed-ledger protocols to see this.

A common element of the first three problems is the demand for \emph{symmetric} protocols, namely, protocols in which all processors run identical local programs. In all three, we will see the power of randomness in action in distributed protocols, which allows them, in different ways, to break this symmetry. We note that such use of randomness was a key element in practical distributed systems, most famously the ALOHA and Ethernet networks in the 1970s~\cite{Abramson1970ALOHA,MetcalfeBoggs1976}. Rabin's papers, together with the work of Francez and Rodeh~\cite{FrancezRodeh1980}, are major initial steps in the \emph{formal} design and analysis of probabilistic distributed protocols in the early 1980s.


\subsection{Dining philosophers}

Edsger Dijkstra, who laid the foundations of distributed computing and concurrent programming, liked to invent stories that abstractly capture the many issues, models, and goals arising in the field.
The \emph{dining philosophers} problem is among the most famous of these stories and has served as a model problem for concurrency control ever since. Dijkstra described it in~\cite{Dijkstra1971Hierarchical}, and Hoare formalized it in~\cite{Hoare1978CSP}. This problem abstracts the difficulties of resource allocation in symmetric environments. We repeat Dijkstra's informal version, which, of course, requires translation into the dry language of processors and resources.

Imagine some philosophers sitting around a table. The philosophers are identical; namely, they have no unique names or IDs. Between each pair of adjacent philosophers is a fork\footnote{For this particular story, some people prefer ``chopstick'' to ``fork,'' but we stick to the original version.}. In front of every philosopher is a bowl containing an unlimited amount of pasta. Bowls are private, but forks are shared. Philosophers alternate between two activities: thinking and eating. To think, no external resource is needed. To eat, a philosopher must use both forks, the one on her left and the one on her right. The only communication mechanism between (adjacent) philosophers is an \emph{atomic} action that tests whether an adjacent fork is available and, if so, takes it. A philosopher who obtains both forks eats for a while and then releases the forks one after the other.

The original challenge calls for a \emph{deadlock-free} program (in which at least one hungry philosopher eventually gets to eat). A stronger requirement that we will also discuss is a \emph{starvation-free} program (in which every hungry philosopher eventually gets to eat). As the setting is symmetric, we demand that the protocol be symmetric as well. Namely, the program for each philosopher is the same. Initially all philosophers start in the same state and all forks are down.

To illustrate the difficulty, consider the following naive program, which philosophers execute whenever they become hungry.

\begin{enumerate}[label=\arabic*.]
\item Try to obtain the left fork.
\item Try to obtain the right fork.
\item If unsuccessful in either attempt, return to step 1.
\item Eat.
\item Release the left fork.
\item Release the right fork.
\end{enumerate}

It is clear that this program is not even deadlock-free. Imagine a scheduler that activates each philosopher for just one step (in any order). After this, each is holding her left fork, and this state will never change, regardless of the future activation order.
Of course, one can try to avoid deadlock, for example, by modifying the program so that it releases the left fork if the right one is not available. A similar scheduler can easily make them all starve again. Can you find a starvation-free program? There were many suggested programs, but all went outside the model above and either broke symmetry (for example, by giving the philosophers different IDs) or changed the problem (for example, by making the philosophers ``hygienic''~\cite{ChandyMisra1984}, or by adding a ``waiter''). But Dijkstra's original problem remained open.

In 1981, Lehmann and Rabin~\cite{LehmannRabin1981} completely resolved this basic problem, both negatively and positively. They proved three main theorems. The first explains past failures, showing that there is no deterministic solution whatsoever.

\begin{theorem}[Lehmann--Rabin~\cite{LehmannRabin1981}]
There is no deadlock-free deterministic, symmetric program for the dining philosophers problem.
\end{theorem}

A proof must present, for \emph{every} possible symmetric program, a scheduler that causes all philosophers to starve. The idea of their proof is simple (though formalizing it is not) and uses the \emph{same} scheduler for all programs. The scheduler behaves as follows: it activates each philosopher in turn for one step, in cyclic order. This constitutes one phase. It then continues, repeating such phases forever. The invariant maintained inductively on the number of phases is that, before each phase begins, the state of the system is symmetric (i.e., all philosophers are in the same state). This is true at the start, and the symmetry of the program ensures that, in each phase, the result of the (inductively identical) step of every philosopher is the same. To conclude, assume for contradiction that in some phase some philosopher is eating. The inductive claim, and the fact that each philosopher executes exactly the same step with the same result in each phase, imply that, at the end of this phase, all philosophers are eating. But this is impossible, as each needs two forks to eat! Observe that asynchrony is crucial to this impossibility result, but it is used here in a rather benign way. The scheduler essentially prevents any philosopher from executing two successive steps of its program before every other philosopher executes one intervening step.


The next two theorems show how to solve the problem, first for deadlock-freedom and then for starvation-freedom, if we allow programs to be probabilistic.

\begin{theorem}[Lehmann--Rabin~\cite{LehmannRabin1981,DuflotFribourgPicaronny2004}]
There is a probabilistic, symmetric program for the dining philosophers problem that is deadlock-free with probability one.
\end{theorem}

Actually, the probabilistic program is a small variation of the naive deterministic program above. Indeed, essentially the only change is that, when a philosopher who holds no forks tries to get one, she tosses a fair coin to decide whether to try the left or the right fork first. The correctness proof is quite nontrivial. The paper makes crucial use of a natural assumption: every philosopher must be activated infinitely often by the scheduler. This assumption was later removed via a modification of the Lehmann--Rabin protocol by Duflot, Fribourg, and Picaronny~\cite{DuflotFribourgPicaronny2004}. They also developed a very elegant formal framework in which to express the analysis of such probabilistic programs.


We conclude with the final gem from~\cite{LehmannRabin1981}. Lehmann and Rabin observe that, for the probabilistic program above, there do exist schedules that prevent all but one of the philosophers from eating. While deadlock-free, this program is far from satisfactory (at least for hungry philosophers). They then devise a new probabilistic program that fixes this problem and satisfies a stronger {\em fairness} condition called starvation-freeness. A key ingredient in the new protocol is a slight variation in the information held by the shared variable represented by a fork. The state of a fork now tells us more than whether or not it is in use. It also records, with a couple of extra bits, which of the two adjacent philosophers used it last, as well as the intention of a philosopher to pick it up. Roughly speaking, the resulting probabilistic program (which they call the ``courteous philosophers' algorithm'') prevents the last eater adjacent to a fork from accessing it next, unless the other is not interested. This modification guarantees that, in every schedule, every philosopher who is activated infinitely often will also get to eat infinitely often with probability one.

\begin{theorem}[Lehmann--Rabin~\cite{LehmannRabin1981}]
There is a probabilistic symmetric program for the dining philosophers problem that is starvation-free with probability one.
\end{theorem}


\subsection{Mutual exclusion}




Another basic resource-allocation problem is mutual exclusion: how to ensure that no more than one of many competing processes can obtain exclusive access to a resource at any given time. The only communication mechanism available to the processes for deciding whether one of them can safely use the resource is a single common \emph{test-and-set} variable: a variable that a single process can \emph{atomically} read and then update in a single step. 

For example, assume the variable takes 2 possible values (namely, 1 bit), ``free'' and ``occupied.'' A possible protocol for the processes is as follows. Each process that wants to gain access to the resource tries to access the variable. If it succeeds and the value is ``free,'' it immediately sets it to ``occupied,'' then uses the resource, and, when done, flips the variable value back to ``free.'' Regardless of which process is allowed access to this variable at each step, it is clear that mutual exclusion is achieved; no two processes can use the resource simultaneously. Also, it achieves a basic fairness condition, namely, that at least one process needing the resource will get access. However, it may be that two processes need the resource, and scheduling of accesses will allow one to repeatedly use it while the other never does.

What we seek is a stronger notion of fairness, called a \emph{bounded-waiting} protocol for this problem. Such a protocol guarantees that no process trying to gain access to that resource waits forever. As before, the environment is asynchronous, with an arbitrary scheduler dictating the timing of the processes' actions and access to the shared variable. In this model, processes may be assumed to have unique IDs which can be used when setting the shared variable. Of course, in this case this variable will need to have more states (or bits). The main parameter to minimize is the number of bits that this shared variable should take, as a function of~$n$, the number of processes.


This model was first studied by Burns, Fischer, Jackson, Lynch, and Peterson~\cite{BurnsFischerJacksonLynchPeterson1982}. They gave a \emph{deterministic} protocol using only $\log n +O(1)$ bits to achieve bounded waiting, and also proved that this is best possible (in other words, the variable must have at least as many states as there are processors). Rabin naturally attempted to beat this bound using randomness. In his paper~\cite{Rabin1982MutualExclusion}, Rabin devised a probabilistic protocol that reduces this bit-length requirement exponentially, to $O(\log\log n)$, guaranteeing bounded waiting with probability one. Later papers~\cite{HartSharirPnueli1983,Saias1992} discovered subtle flaws in Rabin's analysis, which implicitly assumed very weak adversaries. These flaws were corrected, with a modified protocol and new analysis by Kushilevitz and Rabin~\cite{KushilevitzRabin1992}.


\begin{theorem}[Kushilevitz--Rabin~\cite{KushilevitzRabin1992}]
Mutual exclusion with bounded waiting can be achieved with probability one by a probabilistic symmetric protocol using only $O(\log\log n)$ bits.
\end{theorem}

With Kushilevitz, Mansour, and Zuckerman, Rabin proceeded to prove that the above result is best possible~\cite{KushilevitzMansourRabinZuckerman1998}.

\begin{theorem}[Kushilevitz--Mansour--Rabin--Zuckerman~\cite{KushilevitzMansourRabinZuckerman1998}]
Any bounded-waiting probabilistic symmetric protocol for mutual-exclusion must use $\Omega(\log\log n)$ bits.
\end{theorem}


\subsection{Byzantine generals}

Another fundamental problem in distributed systems is Byzantine agreement, which captures \emph{coordination} in the presence of faults. Roughly, a collection of processors must agree on a common value, even though some of them may be faulty. The problem was formulated by Pease, Shostak, and Lamport~\cite{PeaseShostakLamport1980,LamportShostakPease1982}. Their abstraction imagines Byzantine generals in a battle, facing a critical decision of whether to attack or retreat. Each has its own view, but it is paramount that they all agree on a common action. They can send messages (for example, via couriers) to one another to achieve this coordinated decision. However, the problem is complicated by a variety of possible faults. Messages may be arbitrarily delayed. Also, some generals are faulty: they may die or, worse, be corrupt and send inconsistent messages to different parties. Indeed, the word \emph{Byzantine} is meant to capture the worst possible scenario: an adversary may arbitrarily schedule the arrival of messages and may choose the corrupt generals arbitrarily. Of course, we do not expect anything from the faulty generals, but we want all nonfaulty generals eventually to agree on an action in finite time for any adversary. To prevent a trivial solution in which they all choose an action (for example, ``attack'') as a default, regardless of their initial values, we demand that, if they all have the same value, that value must be the chosen action. As usual, we seek a symmetric protocol.


This problem has been studied in both the asynchronous and synchronous models. Two parameters are important here: $n$, the number of generals, and $t$, the maximum number of faulty ones. We summarize the results in both, comparing deterministic and probabilistic protocols.

In the asynchronous model, a celebrated result of Fischer, Lynch, and Paterson~\cite{FischerLynchPaterson1985} shows that, even if there is only $t=1$ fault (and even if this fault is simply a ``crash fault,'' for example, one general stops communicating at some point), as long as $n\geq 2$, there is no deterministic protocol for this problem. Here again, as with dining philosophers, randomization makes the impossible possible. Ben-Or, soon after obtaining his PhD under Rabin, devised a probabilistic symmetric protocol that can tolerate any $t<n/5$ faults and terminates in finite expected time~\cite{BenOr1983}. This was improved to (the optimal!) $t<n/3$ by Bracha~\cite{Bracha1987}. Unfortunately, both algorithms may take $\exp(n)$ time\footnote{Time must be defined in an asynchronous model; here it is taken to be the length of the longest sequence of message deliveries.} to terminate when the number of faults~$t$ is linear in~$n$. Breakthrough work of King and Saia~\cite{KingSaia2013,KingSaia2018Correction} devised a protocol that tolerates $t=\Theta(n)$ faults in polynomial expected time. More recently, Huang, Pettie, and Zhu achieved polynomial expected time at the optimal resilience threshold $t<n/3$~\cite{HuangPettieZhu2024}.


In the synchronous model, there is a global clock for all processors, who proceed in rounds. Messages sent in a round are delivered by the end of that round. Here, deterministic symmetric protocols do exist. However, the number of rounds of communication is essentially the maximum number of faults. The original papers~\cite{PeaseShostakLamport1980,LamportShostakPease1982} prove that $n>3t$ is necessary, and give a $(t+1)$-round protocol for the problem when this condition holds. The tight lower bound of $t+1$ on the number of rounds required was proved by Fischer and Lynch~\cite{FischerLynch1982}. These results are in the information-theoretic model, where the adversary can understand all messages and impersonate any general. Dolev and Strong~\cite{DolevStrong1983} proved that even in the cryptographic setting, where in particular we allow authenticated messages,\footnote{These prevent, for example, $A$ from falsely claiming, ``I received message $m$ from $B$.''}, still $t+1$ rounds are necessary. Clearly, when $t$ is linear in~$n$, this may be impractical.


Rabin shatters this bound using randomization. He devises a symmetric probabilistic protocol that solves the problem in an expected number of rounds that is \emph{constant}, independent of~$t$. A central point (also in Ben-Or's result above) is that even though the Byzantine adversary can completely observe the protocol, it cannot predict future coin tosses!

\begin{theorem}[Rabin~\cite{Rabin1983Byzantine}]
The Byzantine generals problem has a probabilistic protocol that, for every $t<n/10$, terminates in four rounds, in expectation. Moreover, the probability that the number of rounds exceeds any number~$k$ is at most $\exp(-k)$.
\end{theorem}


The above theorem is actually two theorems, about two protocols. When we allow the generals only crash faults, a simple protocol satisfies the theorem in the information-theoretic setting (Rabin calls this variant the \emph{distributed commit} problem). For arbitrary Byzantine faults, a much more involved protocol satisfies the theorem only in the cryptographic setting (using authentication) and, moreover, makes another strong assumption. Rabin assumes the availability of special correlated random strings to each processor before the protocol starts\footnote{Essentially, these are shares, in Shamir's secret-sharing scheme~\cite{Shamir1979}, of random secrets allocated to each processor.}. In his paper~\cite{Rabin1983Byzantine}, Rabin poses the challenge of removing this assumption, namely, of obtaining such correlated randomness organically within the protocol. 

This challenge was met in the important paper by Feldman and Micali~\cite{FeldmanMicali1988}. They use \emph{verifiable secret sharing}, an extension of Shamir's classical secret-sharing scheme~\cite{Shamir1979}, together with a variety of new tools, to give a constant-round protocol to generate shares of such secrets. Thus they remove the assumption in Rabin's paper, satisfying the theorem in the cryptographic setting\footnote{If private links between every pair of players exist, no cryptography is needed.} under Byzantine faults; moreover, they improve his tolerance bound to the optimal $t<n/3$.


\subsection{Information dispersal}

We conclude with another beautiful gem of Rabin. Consider the extremely basic problem, in distributed systems with faults, of transmitting information. Imagine the processors connected through some sparse network~$N$. One processor, $A$, would like to send a file~$F$ of length~$k$ to processor~$B$. If there are no faults, the file can be sent along some $A$--$B$ path in the graph. But in the presence of faults, the failure of any node or edge on that path would cause the information to be lost.

There are many possible fault models, both probabilistic and adversarial. From any such model, we need just one parameter: $m$, an upper bound on the number of disjoint paths that can fail simultaneously during the transmission of one message.


There is a trivial solution that guarantees safe delivery: send~$F$ along $m+1$ disjoint $A$--$B$ paths. Surely, one copy of~$F$ will arrive intact. Of course, the number of bits sent is now $(m+1)k$, a factor of~$m$ larger, which increases the work and congestion in the network by a factor of~$m$. Can this be avoided?

Rabin's ingeniously simple solution~\cite{Rabin1989IDA} enables this. Instead of using $m+1$ paths, he uses $2m$, only a factor of two more. But this enables him to shrink the total amount of information sent to only $2k$, and hence to incur only a factor-two increase in work and congestion, independent of~$m$. Here is how the IDA (Information Dispersal Algorithm) works.

Divide the file~$F$ into $m$ disjoint parts, $F_1,F_2,\ldots,F_m$, each of size $k/m$. We can think of each part as an element of a finite field. Now encode this vector of elements, using a linear \emph{erasure code}, as a vector of $2m$ elements $E_1,E_2,\ldots,E_{2m}$, such that all the $F_i$'s can be recovered from \emph{any} subset of $m$ of the $E_j$'s. (We know many such codes, which also have fast encoding and decoding algorithms; Rabin suggested using a Cauchy matrix, but one can also use a Reed--Solomon code, for example.)

Now everything is obvious. $A$ sends each of the $E_i$'s along one of the $2m$ disjoint $A$--$B$ paths. As at most $m$ of the paths will fail, at least $m$ of the $E_j$'s will arrive at~$B$, and she will be able to reconstruct the file~$F$ from them. Notice that the total amount of information sent is indeed $2k$. Rabin also shows how to implement two variations. One allows the redundancy factor~$2$ to be reduced to any $1+\epsilon$, at the cost of using $m/\epsilon$ paths. The other allows the use of paths with low congestion (namely every edge is used by only a few of them) instead of completely disjoint paths.

This brings us to an issue we have not yet discussed: do such paths exist, and how does one find them in a distributed manner? This of course depends on the network. Rabin suggests using the Boolean hypercube as the underlying network, and the then-recent \emph{randomized routing} scheme of Valiant~\cite{Valiant1982Routing} to find many low-congestion paths\footnote{Another beautiful example of the power of randomness.}. As these are low-congestion rather than disjoint, Rabin needs to explain how to implement the IDA in this setting. By now, this problem is extremely well understood; for a variety of highly \emph{expanding} sparse networks there are even deterministic distributed protocols for finding such paths.


\vspace{.5in}
\phantomsection
\addcontentsline{toc}{section}{Acknowledgements}

\noindent
{\bf \large Acknowledgements}

Many aided me in writing this survey, and I would like to acknowledge and thank them.

On the computer side, I was helped in many ways by ChatGPT 6. I used it in the obvious ways to locate, organize, and summarize material so that I could focus on what I needed to read. I used it to correct typos and grammar, fix \LaTeX{} errors, and create the bibliography files. I used it to go over my drafts and suggest revisions on historical, technical, and other grounds. But we also had many fun conversations, and arguments, about the best presentation of the material. 

On the human side, many read all or parts of the text and gave me excellent comments that improved it: Michael Ben-Or, Nataly Brukhim, Irit Dinur, David Harel, Eyal Kushilevitz, Nati Linial, Salil Vadhan, Les Valiant, Edna Wigderson, Yuval Wigderson. Thank you all!

\clearpage
\phantomsection
\addcontentsline{toc}{section}{References}
\bibliographystyle{bib/surveyplain}
\bibliography{bib/survey}

\clearpage
\appendix
\section{Complete list of Rabin's publications and technical reports}
\label{app:rabin-publications}

The following list is based on Rabin's Harvard publication list, supplemented and checked where possible against bibliographic databases, publisher records, and institutional archives. Distinct technical-report, conference, and journal versions are retained when they circulated separately or appeared in different years. Entries are ordered by publication year, and alphabetically by title within each year. Each separately listed technical-report, conference, or journal version is placed according to its own year of publication or circulation.

\begingroup
\small
\setlength{\parskip}{0pt}
\begin{enumerate}[leftmargin=*,label=\arabic*.,itemsep=0.45em,parsep=0pt,topsep=0.4em]
\item Michael O. Rabin. \emph{A theorem on partially ordered sets (Hebrew)}. Riveon Lematematika, vol. 7 (1953), pp. 26-29.
\item Michael O. Rabin. \emph{Sur la repr{\'e}sentation des id{\'e}aux par des id{\'e}aux primaires}. C. R. de l'Acad{\'e}mie des sciences, vol. 237 (1953), pp. 544-545.
\item Michael O. Rabin. \emph{On regular polygons with lattice-point vertices (Hebrew)}. Riveon Lematematika, vol. 8 (1954), pp. 13-15.
\item Michael O. Rabin. \emph{A note on Helly's theorem}. Pacific Journal of Mathematics, vol. 5 (1955), pp. 363-366.
\item Michael O. Rabin. \emph{Effective computability of winning strategies}. Annals of Mathematics Studies, vol. 39 (1957), pp. 147-157.
\item Michael O. Rabin. \emph{Two way finite automata}. Proceedings of the Summer Institute of Symbolic Logic, 1957 at Cornell, pp. 366-369.
\item Michael O. Rabin. \emph{Recursive unsolvability of group theoretic problems}. Annals of Mathematics, vol. 67 (1958), pp. 172-194.
\item Michael O. Rabin and Robert Z. Norman. \emph{An algorithm for a minimum cover of a graph}. Proceedings of the American Mathematical Society, vol. 10 (1959), pp. 315-319.
\item Michael O. Rabin. \emph{Arithmetical extensions with prescribed cardinality}. Indagationes Mathematicae, vol. 21 (1959), pp. 439-446.
\item Michael O. Rabin and Dana S. Scott. \emph{Finite automata and their decision problems}. IBM Journal, vol. 3 (1959), pp. 114-125. Russian translation: Konechnye avtomaty i zadachi ikh razresheniya, Kiberneticheskii sbornik, vol. 4 (1962), pp. 58-91. Reprinted in Sequential Machines, Selected Papers, E. F. Moore (editor), Addison-Wesley, Reading, MA (1964), pp. 63-92.
\item Michael O. Rabin and W. Wesley Peterson. \emph{On codes for checking logical operations}. IBM Journal, vol. 3 (1959), pp. 163-168.
\item Michael O. Rabin. \emph{On recursively enumerable and arithmetic models of set theory}. Journal of Symbolic Logic, vol. 23 (1958), appeared 1959, pp. 408-416.
\item Michael O. Rabin. \emph{Speed of computation and classification of recursive sets}. Third Convention of Scientific Societies, Israel (1959), pp. 1-2.
\item Michael O. Rabin. \emph{Computable algebra, general theory and theory of computable fields}. Transactions of the American Mathematical Society, vol. 94 (1960), pp. 341-360.
\item Michael O. Rabin. \emph{Degree of difficulty of computing a function and a partial ordering of recursive sets}. Technical Report O.N.R. Contract (1960), pp. 341-360.
\item Michael O. Rabin. \emph{Non-standard models and the independence of the induction axiom}. Essays on the Foundations of Mathematics, dedicated to A. A. Fraenkel (1961), pp. 287-299.
\item Michael O. Rabin. \emph{Classes of structures and sets of sentences with the intersection property}. Actes du Colloque de mathematiques a l'Occasion du Tricentenaire de la Mort de B. Pascal, Tome 1 (1962), pp. 39-53.
\item Michael O. Rabin. \emph{Diophantine equations and non-standard models of arithmetic}. Proceedings of the International Congress for Logic and Methodology of Science, Stanford (1962), pp. 151-158.
\item Michael O. Rabin. \emph{Probabilistic Automata}. Information and Control, vol. 6 (1963), pp. 230-245. Reprinted in Sequential Machines, Selected Papers, E. F. Moore (editor), Addison-Wesley, Reading, Mass., pp. 98-114.
\item Michael O. Rabin. \emph{Real time computation}. Israel Journal of Mathematics, vol. 1 (1963), pp. 203-211.
\item Micha A. Perles, Michael O. Rabin, and Eli Shamir. \emph{The theory of definite automata}. IEEE Transactions on Electronic Computers, vol. EC-12 (1963), pp. 233-243.
\item Michael O. Rabin and Hao Wang. \emph{Words in the history of a Turing machine with a fixed input}. Journal of the ACM, vol. 10 (1963), pp. 526-527.
\item Michael O. Rabin. \emph{A simple method for undecidability proofs}. Proceedings of the 1964 International Congress for Logic, North-Holland (1965), pp. 58-68.
\item Michael O. Rabin. \emph{Universal groups of automorphisms of models}. Berkeley International Symposium on the Theory of Models, North-Holland (1965), pp. 274-284.
\item Michael O. Rabin. \emph{Classical and probabilistic automata}. Automata Theory, E. B. Cainiello (editor), Academic Press, 1966, pp. 304-313.
\item Calvin C. Elgot and Michael O. Rabin. \emph{Decidability and undecidability of extensions of second (first) order theory of (generalized) successor}. Journal of Symbolic Logic, vol. 31 (1966), pp. 169-181.
\item Michael O. Rabin. \emph{Decidability of Second Order Theories and Automata on Infinite Trees}. Bulletin of the American Mathematical Society, vol. 74 (1968), pp. 1025-1029.
\item Michael O. Rabin. \emph{Mathematical theory of automata}. Proceedings of Symposia in Applied Mathematics, vol. 19, American Mathematical Society, Providence, RI (1968), pp. 153-175.
\item Michael O. Rabin. \emph{Decidability of second-order theories and automata on infinite trees}. Transactions of the American Mathematical Society, vol. 141 (1969), pp. 1-35.
\item Michael O. Rabin. \emph{Automata on infinite trees and the synthesis problem}. ONR Technical Report 37, Jerusalem (1970), 24 pp.
\item Michael O. Rabin. \emph{Weakly definable relations and special automata}. Mathematical Logic and Foundations of Set Theory, Y. Bar-Hillel (editor), North-Holland (1970), pp. 1-23.
\item Michael O. Rabin. \emph{Decidability and definability in second-order theories}. Proceedings of the International Mathematical Congress (1970), published 1971.
\item Michael O. Rabin. \emph{Automata on infinite objects and Church's problem}. Regional Conference Series in Mathematics, No. 13, American Mathematical Society (1972), 22 pp.
\item Michael O. Rabin and Shmuel Winograd. \emph{Fast evaluation of polynomials by rational preparation}. Communications on Pure and Applied Mathematics, vol. 25 (1972), pp. 453-458.
\item Michael O. Rabin. \emph{Proving simultaneous positivity of linear forms}. Journal of Computer and System Sciences, vol. 6 (1972).
\item Michael O. Rabin. \emph{Solving linear equations by means of scalar products}. Complexity of Computer Computations, R. E. Miller (editor), Plenum Press (1972), pp. 11-20.
\item Vaughan R. Pratt, Michael O. Rabin and Larry J. Stockmeyer. \emph{A Characterization of the Power of Vector Machines}. Proceedings of the Sixth Annual ACM Symposium on Theory of Computing (STOC), ACM (1974), pp. 122-134. DOI: \href{https://doi.org/10.1145/800119.803892}{\nolinkurl{10.1145/800119.803892}}.
\item Michael J. Fischer and Michael O. Rabin. \emph{Super exponential complexity of Presburger arithmetic}. In Complexity of Computations, SIAM-AMS Proceedings, vol. 7 (1974), pp. 27-41.
\item Michael O. Rabin. \emph{Theoretical impediments to artificial intelligence}. Proceedings of the 1974 IFIP Congress, North- Holland, pp. 615-619.
\item Michael O. Rabin. \emph{Probabilistic algorithms}. In Algorithms and Complexity: New Directions and Recent Trends, J. F. Traub (editor), Academic Press, New York (1976), pp. 21-39.
\item Michael O. Rabin. \emph{Complexity of computations}. Communications of the ACM, vol. 20 (1977), pp. 625-633. DOI: 10.1145/359810.359816.
\item Michael O. Rabin. \emph{Decidable theories}. Handbook of Mathematical Logic, J. Barwise (editor), North-Holland, Amsterdam (1977), pp. 595-629.
\item Michael O. Rabin. \emph{Digitalized signatures}. Foundations of Secure Computations, R. De Millo and R. Lipton (editors), Academic Press, New York (1978), pp. 155-166.
\item Michael O. Rabin. \emph{Digitalized signatures and public-key functions as intractable as factorization}. MIT/LCS/\allowbreak TR-212, 1979.
\item Walter Baur and Michael O. Rabin. \emph{Linear Disjointness and Algebraic Complexity}. L'Enseignement Mathematique, vol. 26 (1980), pp. 333-344.
\item Michael O. Rabin. \emph{Probabilistic algorithms in finite fields}. SIAM Journal on Computing, vol. 9 (1980), pp. 273-280.
\item Michael O. Rabin. \emph{Probabilistic tests for primality}. Journal of Number Theory, vol. 12 (1980), pp. 128-138.
\item Daniel Lehmann and Michael O. Rabin. \emph{A Symmetric and Fully Distributed Solution to the Dining Philosophers Problem}. (Extended Abstract), 8th ACM Symposium on Principles of Programming Languages (1981), pp. 133-138.
\item Richard M. Karp and Michael O. Rabin. \emph{Efficient Randomized Pattern-Matching Algorithms}. Aiken Computation Laboratory, Harvard University, TR-31-81, Cambridge, MA, 1981. Technical-report version; later journal version appears as RabinPub054. Cited in later systems/security bibliographies as Aiken Laboratory TR-31-81.
\item Michael O. Rabin. \emph{Fingerprinting by random polynomials}. Center for Research in Computing Technology, Harvard University, TR-15-81, 1981.
\item Michael O. Rabin. \emph{How to Exchange Secrets with Oblivious Transfer}. Technical Report TR-81, Aiken Computation Laboratory, Harvard University, 1981; also Cryptology ePrint Archive, Paper 2005/187, 2005.
\item Michael O. Rabin. \emph{N-process mutual exclusion with bounded waiting by 4 log2 N-valued shared variable}. Journal of Computer and System Sciences, vol. 25 (1982), pp. 66-75.
\item Michael O. Rabin. \emph{The choice coordination problem}. Acta Informatica, vol. 17 (1982), pp. 121-134.
\item Michael O. Rabin. \emph{Randomized Byzantine Generals}. IEEE 24th Symposium on Foundations of Computer Science (1983), pp. 403-409.
\item Michael O. Rabin. \emph{Transaction protection by beacons}. Journal of Computer and System Sciences, vol. 27 (1983).
\item Michael O. Rabin. \emph{Discovering repetitions in strings}. Combinatorial Algorithms on Words, Galil and Apostolico (editors), Springer-Verlag, Berlin CS series (1985), pp. 279-289.
\item Michael O. Rabin and Jeffrey Shallit. \emph{Randomized algorithms in number theory}. Communications on Pure and Applied Mathematics, vol. 39 (1986), pp. 239-256.
\item Joseph Y. Halpern and Michael O. Rabin. \emph{A logic to reason about likelihood}. Artificial Intelligence, vol. 32 (1987), pp. 379-405.
\item Benny Chor and Michael O. Rabin. \emph{Achieving Independence in Logarithmic Number of Rounds}. Proceedings of the 6th ACM Conference on Principles of Distributed Computing, August 1987, pp. 260-268.
\item Michael O. Rabin and J. D. Tygar. \emph{An integrated toolkit for operating system security}. Center for Research in Computing Technology, Harvard University, Technical Report TR-05-87, 1987.
\item Michael O. Rabin. \emph{Efficient Dispersal of Information for Security and Fault Tolerance}. Aiken Computation Laboratory, Harvard University, TR-02-87, Cambridge, MA, 1987. Technical-report version of the information-dispersal work; journal and book-chapter versions appear as RabinPub056 and RabinPub060.
\item Richard M. Karp and Michael O. Rabin. \emph{Efficient randomized pattern-matching algorithms}. IBM Journal of Research and Development, vol. 31 (1987), pp. 249-260.
\item Michael O. Rabin and J. D. Tygar. \emph{An Integrated Toolkit for Operating System Security}. Center for Research in Computing Technology, Aiken Computation Laboratory, Harvard University, TR-05-87R, Cambridge, MA, 1988. Revised technical-report version cited in later bibliographies; related to RabinPub091 and RabinPub059.
\item Michael O. Rabin and J. D. Tygar. \emph{An Integrated Toolkit for Operating System Security}. In Foundations of Data Organization and Algorithms, W. Litwin and H. J. Schek (editors), 3rd International Conference, FODO 1989, Paris, France, June 1989, Springer-Verlag, pp. 2-15.
\item Noga Alon and Michael O. Rabin. \emph{Biased Coins and Randomized Algorithms}. In Randomness and Computing, S. Micali (editor), Advances in Computing, vol. 5, JAI Press (1989), pp. 499-507.
\item Michael O. Rabin. \emph{Efficient Dispersal of Information for Security, Load Balancing, and Fault Tolerance}. Journal of the ACM, vol. 38 (1989), pp. 335-348.
\item Michael O. Rabin and Vijay V. Vazirani. \emph{Maximum Matchings in General Graphs Through Randomization}. Journal of Algorithms, vol. 10 (December 1989).
\item Michael O. Rabin. \emph{The Information Dispersal Algorithm and Its Applications}. In Sequences: Combinatorics, Compression, Security, and Transmission, Capocelli (editor), Springer-Verlag, New York, 1990, pp. 406-419.
\item Noga Alon, Daniel J. Kleitman, Richard J. Lipton, Roy Meshulam, Michael O. Rabin, and Joel Spencer. \emph{Set Systems with No Union of Cardinality 0 Modulo m}. Graphs and Combinatorics, vol. 7 (1991), pp. 97-99.
\item Michael O. Rabin. \emph{Dependable Parallel Computing by Randomization (Abstract)}. Future Tendencies in Computer Science, Control and Applied Mathematics: International Conference on the Occasion of the 25th Anniversary of INRIA, Lecture Notes in Computer Science, vol. 653, Springer (1992), pp. 17-20. DOI: \href{https://doi.org/10.1007/3-540-56320-2_50}{\nolinkurl{10.1007/3-540-56320-2_50}}.
\item Zvi M. Kedem, Krishna V. Palem, Michael O. Rabin, and A. Raghunathan. \emph{Efficient Program Transformations for Resilient Parallel Computation via Randomization}. Efficient Program Transformation for Resilient Computation via Randomization, Proceedings of the Annual ACM Symposium on Theory of Computing (STOC), 1992, pp. 306-318.
\item Yonatan Aumann and Michael O. Rabin. \emph{Fast PRAM Simulation on Fully Asynchronous Parallel Systems}. Proceedings of the 33rd Annual Symposium on Foundations of Computer Science (FOCS), October 1992, pp. 147-156.
\item Eyal Kushilevitz and Michael O. Rabin. \emph{Randomized Mutual Exclusion Algorithms Revisited}. Proceedings of Principles of Distributed Computing (PODC), August 1992, pp. 235-283.
\item Yonatan Aumann, Zvi M. Kedem, Krishna V. Palem, and Michael O. Rabin. \emph{Highly efficient asynchronous execution of large-grained parallel programs}. Proceedings of the 34th IEEE Annual Symposium on Foundations of Computer Science (FOCS), 1993, pp. 271-280.
\item Eyal Kushilevitz, Yishay Mansour, Michael O. Rabin, and David Zuckerman. \emph{Lower Bounds for Randomized Mutual Exclusion (Extended Abstract)}. Proceedings of the 25th Annual Symposium on Theory of Computing (STOC), May 16-18, 1993, pp. 154-163.
\item Michael O. Rabin. \emph{Optimal Parallel Pattern Matching Through Randomization}. Sequences II: Methods in Communication, Security, and Computer Science, R. Capocelli et al. (editors), Springer-Verlag, 1993, pp. 292-299.
\item Yonatan Aumann and Michael O. Rabin. \emph{Clock Construction in Fully Asynchronous Parallel Systems and PRAM Simulation}. Theoretical Computer Science, vol. 128 (1994), pp. 3-30.
\item Daniel Lehmann and Michael O. Rabin. \emph{The Advantages of Free Choice: A symmetric and Fully Distributed Solution for the Dining Philosophers Problem}. In A Classical Mind: Essays in Honour of C. A. R. Hoare, A. W. Roscoe (editor), Prentice-Hall International Series in Computer Science (1994), pp. 337-356.
\item Eyal Kushilevitz, Yishay Mansour, and Michael O. Rabin. \emph{On Lotteries with Unique Winners}. SIAM Journal on Discrete Mathematics, vol. 8, no. 1, pp. 93-98, February 1995.
\item Partha Dasgupta, Zvi M. Kedem, and Michael O. Rabin. \emph{Parallel Processing on Networks of Workstations: A Fault-Tolerant, High Performance Approach}. Proceedings of the 15th International Conference on Distributed Computing Systems, 1995 (IEEE Outstanding Paper Award).
\item Michael O. Rabin. \emph{Computationally Hard Algebraic Problems (Invited Lecture)}. Proceedings of the 37th IEEE Annual Symposium on Foundations of Computer Science (FOCS), 1996.
\item Michael O. Rabin. \emph{Correctness of Programs and Protocols through Randomization}. Advances in Computing Science--ASIAN '97, Lecture Notes in Computer Science, vol. 1345, Springer (1997), pp. 1-3. DOI: \href{https://doi.org/10.1007/3-540-63875-X_38}{\nolinkurl{10.1007/3-540-63875-X_38}}.
\item Michael O. Rabin and Silvio Micali. \emph{Hashing on Strings, Cryptography, and Protection of Privacy}. Proceedings Compression and Complexity of Sequences, IEEE Computer Society, Los Alamitos, CA, June 11-13, 1997, p. 1.
\item Yonatan Aumann and Michael O. Rabin. \emph{Authentication, Enhanced Security and Error Correcting Codes, IACR Distinguished Lecture}. In Advances in Cryptology--CRYPTO 98, Lecture Notes in Computer Science 1462, Springer, 1998, pp. 299-303.
\item Rosario Gennaro, Michael O. Rabin, and Tal Rabin. \emph{Simplified VSS and Fast-track Multiparty Computations with Applications to Threshold Cryptography}. Proceedings of the 17th ACM Symposium on Principles of Distributed Computing, ACM, 1998, pp. 101-112.
\item Richard J. Lipton, Laura F. Landweber, and Michael O. Rabin. \emph{DNA2DNA Computations: A Potential Killer Application}. In DNA Based Computers III, DIMACS Series in Discrete Mathematics and Theoretical Computer Science, vol. 48, H. Rubin (editor), American Mathematical Society, 1999, pp. 161-172.
\item Yonatan Aumann and Michael O. Rabin. \emph{Information Theoretically Secure Communication in the Limited Storage Model}. In Advances in Cryptology --CRYPTO 99, Lecture Notes in Computer Science 1666, Springer, 1999, pp. 65-79.
\item Silvio Micali, Michael O. Rabin, and Salil P. Vadhan. \emph{Verifiable Random Functions}. Proceedings of the 40th Symposium on Foundations of Computer Science (FOCS), IEEE Computer Society, 1999, pp. 120-130.
\item Yonatan Aumann, Johan H{\aa}stad, Michael O. Rabin, and Madhu Sudan. \emph{Linear consistency testing}. Journal of Computer and System Sciences, vol. 62, pp. 589-607, 2001.
\item Yonatan Aumann, Yan Zong Ding, and Michael O. Rabin. \emph{Everlasting Security in the Bounded Storage Model}. IEEE Transactions on Information Theory, vol. 48, issue 6, pp. 1668-1680, June 2002.
\item Yan Zong Ding and Michael O. Rabin. \emph{Hyper-Encryption and Everlasting Security}. Proceedings of the 19th International Symposium on Theoretical Aspects of Computer Science (STACS), pages 1-26, Antibes--Juan les Pins, France, March 14-16, 2002.
\item Michael A. Bender and Michael O. Rabin. \emph{Online Scheduling of Parallel Programs on Heterogeneous Systems with Applications to Cilk}. Theory of Computing Systems, special issue on SPAA '00, vol. 35, pp. 289-304, 2002.
\item Michael O. Rabin. \emph{Hyper Encryption and Everlasting Secrets}. Algorithms and Complexity: Fifth Italian Conference (CIAC), Lecture Notes in Computer Science, vol. 2653, Springer (2003), pp. 7-10. DOI: \href{https://doi.org/10.1007/3-540-44849-7_7}{\nolinkurl{10.1007/3-540-44849-7_7}}.
\item Joe Kilian, Silvio Micali, and Michael O. Rabin. \emph{Zero Knowledge Sets}. Symposium on Foundations of Computer Science (FOCS), pp. 80-86, October 2003.
\item Jonathan Katz, Rafail Ostrovsky, and Michael O. Rabin. \emph{Identity-Based Zero Knowledge}. In Security in Communication Networks, 4th International Conference, 2004, Lecture Notes in Computer Science 3352, Springer, 2005, pp. 180-192.
\item Michael O. Rabin and Dennis E. Shasha. \emph{Preventing Piracy While Preserving Privacy}. Dr. Dobb's Journal Computer Security, CMP Media LLC, 2005.
\item Michael O. Rabin. \emph{Provably Unbreakable Hyper-Encryption in the Limited Access Model}. Proceedings IEEE Symposium on Information Theory Workshop on Theory and Practice in Information-Theoretic Security, pp. 34-37, Awaji Island, Japan, 2005.
\item David C. Parkes, Michael O. Rabin, Stuart M. Shieber, and Christopher Thorpe. \emph{Practical Secrecy-Preserving, Verifiably Correct and Trustworthy Auctions}. ICEC '06: Proceedings of the 8th International Conference on Electronic Commerce, pp. 70-81, New York, NY, USA, ACM Press, 2006.
\item Michael O. Rabin. \emph{DISC 20th Anniversary: Invited Talk Provably Unbreakable Hyper-Encryption Using Distributed Systems}. Distributed Computing: 21st International Symposium (DISC), Lecture Notes in Computer Science, vol. 4731, Springer (2007), pp. 506-508. DOI: \href{https://doi.org/10.1007/978-3-540-75142-7_48}{\nolinkurl{10.1007/978-3-540-75142-7_48}}.
\item Michael O. Rabin, Rocco A. Servedio, and Christopher Thorpe. \emph{Highly Efficient Secrecy-Preserving Proofs of Correctness of Computations and Applications}. LICS 2007, Proceedings of the 22nd Annual IEEE Symposium on Logic in Computer Science, 14 pages, 2007.
\item David C. Parkes, Michael O. Rabin, Stuart M. Shieber, and Christopher Thorpe. \emph{Practical Secrecy-Preserving, Verifiably Correct and Trustworthy Auctions}. Electronic Commerce Research and Applications, vol. 7, no. 3, pp. 294-312, 2008.
\item David C. Parkes, Michael O. Rabin and Christopher Thorpe. \emph{Cryptographic Combinatorial Clock-Proxy Auctions}. Financial Cryptography and Data Security: 13th International Conference (FC), Lecture Notes in Computer Science, vol. 5628, Springer (2009), pp. 305-324. DOI: \href{https://doi.org/10.1007/978-3-642-03549-4_19}{\nolinkurl{10.1007/978-3-642-03549-4_19}}.
\item Michael O. Rabin. \emph{Never too early to begin: computer science for high-school students}. Proceedings of the 17th ACM Annual Conference on Innovation and Technology in Computer Science Education (ITiCSE), ACM (2012), p. 1. DOI: \href{https://doi.org/10.1145/2325296.2325298}{\nolinkurl{10.1145/2325296.2325298}}.
\item Yishay Mansour, S. Muthukrishnan, Michael O. Rabin, and Moti Yung. \emph{Strictly Black-Box Zero-Knowledge and Efficient Validation of Financial Transactions}. ICALP (1) 2012, pp. 738-749.
\item Silvio Micali and Michael O. Rabin. \emph{Cryptography Miracles, Secure Auctions, Matching Problem Verification}. Communications of the ACM, vol. 57, no. 2, pp. 85-93, 2014.
\end{enumerate}
\endgroup

\end{document}